\documentclass[11pt,draft]{amsart}
\usepackage{a4,amsmath,amssymb,amscd}
\linethickness{0.5pt}

\newtheorem{prop}[equation]{Proposition}
\newtheorem{thm}[equation]{Theorem}
\newtheorem{cor}[equation]{Corollary}
\newtheorem{lem}[equation]{Lemma}

\theoremstyle{definition}
\newtheorem{exe}[equation]{Exercise}
\newtheorem{exa}[equation]{Example}
\newtheorem{prob}[equation]{Problem}

\theoremstyle{remark}
\newtheorem*{rem}{Remark}

\numberwithin{equation}{section}

\newcommand{\D}{\Delta}
\renewcommand{\L}{\Lambda}

\renewcommand{\b}{\beta}
\newcommand{\e}{\varepsilon}
\renewcommand{\t}{\theta}
\renewcommand{\l}{\lambda}
\newcommand{\s}{\sigma}
\renewcommand{\phi}{\varphi}

\newcommand{\mb}[1]{{\textbf {\textit#1}}}
\newcommand{\sbr}[2]{{\textstyle\genfrac{[}{]}{}{}{#1}{#2}}}
\newcommand{\bin}[2]{{\textstyle\binom{#1}{#2}}}

\newcommand{\cone}{\mathop{\rm cone}}
\newcommand{\link}{\mathop{\rm link}\nolimits}
\renewcommand{\star}{\mathop{\rm star}\nolimits}

\newcommand{\bideg}{\mathop{\rm bideg}}

\newcommand{\hd}{\mathop{\rm hd}}
\newcommand{\Ker}{\mathop{\rm Ker}}
\newcommand{\Hom}{\operatorname{Hom}}
\newcommand{\Tor}{\operatorname{Tor}}
\newcommand{\djs}{\mathop{\mbox{\it DJ\/}}\nolimits}
\newcommand{\zk}{\mathcal Z_K}
\newcommand{\cc}{\mathop{\rm cc}}
\renewcommand{\k}{\mathbf k}
\newcommand{\Z}{\mathbb Z}
\newcommand{\R}{\mathbb R}
\newcommand{\C}{\mathbb C}

\newcommand{\<}{\langle}
\renewcommand{\>}{\rangle}
\renewcommand{\le}{\leqslant}
\renewcommand{\ge}{\geqslant}

\begin{document}

\title{Stanley--Reisner rings and torus actions}
\author{Taras Panov}
\thanks{The author is grateful to Mikiya Masuda of
Osaka City University and Nat\`alia Castellana of
Universitat Aut\`onoma de Barcelona for the insight gained from
numerous discussions and great hospitality. The work was partially supported
by the Russian Foundation for Basic Research, grant no.~01-01-00546.
Preprint id. ITEP-TH-91/02.}
\address{Department of Mathematics and Mechanics, Moscow State University,
Moscow 119992, Russia\newline
\emph{and }Institute for Theoretical and Experimental Physics, Moscow 117259, Russia}
\email{tpanov@mech.math.msu.su}

\begin{abstract}
We review a class of problems on the borders of topology of torus actions,
commutative homological algebra and combinatorial geometry, which is
currently being investigated by Victor Buchstaber and the author. The text
builds on the lectures delivered on the transformation group courses in
Osaka City University and Universitat Aut\`onoma de Barcelona. We start
with discussing several well-known results and problems on combinatorial
geometry of polytopes and simplicial complexes, and then move gradually
towards investigating the combinatorial structures associated with spaces
acted on by the torus. Parallelly, we set up the required commutative
algebra apparatus, including
Co\-hen--\-Ma\-cau\-lay/\-Go\-ren\-stein rings and Stanley--Reisner face
rings of simplicial complexes.  \end{abstract}

\maketitle

\section{Combinatorial motivations: $f$-vectors}
Here we review some well-known combinatorial results and problems concerning
the number of faces of simplicial polytopes and complexes.

Let $K^{n-1}$ be an arbitrary $(n-1)$-dimensional simplicial
complex on the vertex set $[m]=\{1,\dots,m\}$. The
\emph{$f$-vector} of $K$ is defined as $\mb
f(K)=(f_0,\dots,f_{n-1})$ where $f_i$ is the number of
$i$-dimensional simplices of $K$. We also put $f_{-1}(K)=1$. The equivalent
information is contained in the \emph{$h$-vector} $\mb
h(K)=(h_0,\dots,h_n)$, defined from the equation
\begin{equation}
\label{hvector}
  h_0t^n+\ldots+h_{n-1}t+h_n=(t-1)^n+f_0(t-1)^{n-1}+\ldots+f_{n-1}.
\end{equation}

Convex geometry provides an important class of simplicial complexes, the
boundaries of \emph{simplicial polytopes}. Given a simplicial
polytope $Q$, we define its $f$-and $h$-vector as those of its boundary
complex. The famous theorem of Billera--Lee and Stanley (proved in 1980)
gives a complete characterisation of $f$-vectors of simplicial polytopes.

For any two positive integers $a$, $i$ there exists a unique {\it binomial
$i$-expansion\/} of $a$ of the form
$$
  a=\bin{a_i}i+\bin{a_{i-1}}{i-1}+\dots+\bin{a_j}j,
$$
where $a_i>a_{i-1}>\dots>a_j\ge j\ge1$. Define
$$
  a^{\langle i\rangle}=
  \bin{a_i+1}{i+1}+\bin{a_{i-1}+1}i+\dots+\bin{a_j+1}{j+1},\quad
  0^{\langle i\rangle}=0.
$$

\begin{exa}
For any $a>0$ we have $a^{\langle 1\rangle}=\binom{a+1}2$.
Let $a=28$, $i=4$. The corresponding binomial expansion is
$$
  28=\bin 64+\bin 53+\bin 32.
$$
Hence,
$$
  28^{\langle 4\rangle}=\bin 75+\bin 64+\bin 43=40.
$$
\end{exa}

\begin{thm}[Billera--Lee, Stanley]
\label{gth}
An integer sequence $(f_0,f_1,\dots,f_{n-1})$ is the $f$-\-vec\-tor of a
simplicial $n$-polytope if and only if the corresponding sequence
$(h_0,\ldots,h_n)$
determined by~{\rm(\ref{hvector})} satisfies the following three conditions:
\begin{enumerate}
\item[(a)] $h_i=h_{n-i}$, \ $i=0,\ldots,n$ (the Dehn--Sommerville equations);
\item[(b)] $h_0\le h_1\le\dots\le h_{\sbr n2}$;
\item[(c)] $h_0=1$, $h_{i+1}-h_i\le(h_i-h_{i-1})^{\langle i\rangle}$, \
$i=1,\ldots,\sbr n2-1$.
\end{enumerate}
\end{thm}

\begin{exa}\label{ldimgth}
1. The first inequality $h_0\le h_1$ from part~(b) of $g$-theorem is
equivalent to $f_0=m\ge n+1$. This just expresses the fact that it takes at
least $n+1$ hyperplanes to bound a polytope in~$\R^n$.

2. The first inequality from part~(c) is equivalent to the upper bound
$$
  f_1\le\bin{f_0}2,
$$
which says that any two vertices are joined by at most one edge.

3. The second inequality $h_1\le h_2$ (for $n\ge 4$) from part~(b)
is equivalent to the lower bound
$$
  f_1\ge nf_0-\bin{n+1}2.
$$
\end{exa}

\begin{exe}
Deduce the following ``$f$-vector form" of the Dehn--\-Som\-mer\-vil\-le
equations
\begin{equation}
\label{DSf}
  f_{k-1}=\sum_{j=k}^n(-1)^{n-j}\bin jk f_{j-1},\quad k=0,1,\ldots,n.
\end{equation}
\end{exe}

An integral sequence $(k_0,k_1,\ldots,k_r)$ satisfying $k_0=1$
and $0\le k_{i+1}\le k_i^{\langle i\rangle}$ for $i=1,\ldots,r-1$ is called
an \emph{$M$-vector} (after M.~Macaulay). One may introduce the
\emph{$g$-vector} $\mb g(K^{n-1})=(g_0,g_1,\ldots,g_{\sbr n2})$ of simplicial
complex $K^{n-1}$ by $g_0=1$, \ $g_i=h_i-h_{i-1}$, \
$i>0$. Then conditions~(b) and~(c) from the $g$-theorem
are equivalent to that the
$g$-vector of a simplicial $n$-polytope is an
$M$-vector. The notion of $M$-vector arises in the
following theorem of commutative algebra.

\begin{thm}[Macaulay]
\label{mvect}
An integral sequence $(k_0,k_1,\ldots,k_r)$ is an $M$-\-vec\-tor if and only if
$$
  k_i=\dim_\k A^{2i},\quad i=1,\ldots,r.
$$
for some commutative graded connected $\k$-algebra $A=A^0\oplus
A^2\oplus\dots\oplus A^{2r}$ generated by degree-two elements.
\end{thm}

We sketch the proof of the necessity part of Theorem~\ref{gth} in the end of
section~4.1. Below we review two more famous results concerning the number of
faces, \emph{Upper Bound} and \emph{Lower Bound theorems} for simplicial
polytopes. Both can be deduced from Theorem~\ref{gth}.

A simplicial polytope $Q$ is called \emph{$k$-neighbourly} if any
$k$ vertices span a face.

\begin{exe}
Using the Dehn--Sommerville equations, show that if $Q$ is a $k$-neighbourly
simplicial $n$-polytope and $k>\sbr n2$, then $Q$ is an
$n$-simplex.
\end{exe}

A $\sbr n2$-neighbourly simplicial $n$-polytope is
called \emph{neighbourly}. An example of a neighbourly $n$-polytope
with arbitrary number of vertices is provided by the \emph{cyclic polytope}
$C^n(m)$ defined as the convex hull of any $m$ distinct points $\mb x(t_i)$, \
$t_1<t_2<\ldots<t_m$, on the \emph{moment curve}
$$
  \mb x\colon\R\longrightarrow\R^n,\qquad t\mapsto\mb x(t)=(t,t^2,\ldots,t^n)
  \in\R^n.
$$

\begin{exe}
Prove that the $C^n(m)$ is indeed a neighbourly simplicial
polytope, and that its combinatorial type is independent on a
choice of points on the moment curve.
\end{exe}

\begin{thm}[UBT for simplicial polytopes]
\label{ubt}
  The number of $i$-faces of arbitrary simplicial $n$-polytope $Q$
  does not exceed the number of $i$-faces of any neighbourly $n$-polytope
  with the same number of vertices. That is, if $f_0(Q)=m$, then
  $$
    f_i(Q)\le f_i\bigl(C^n(m)\bigr) \quad\text{ for }i=1,\ldots,n-1.
  $$
  The equality above holds if and only if $Q$ is a
  neighbourly polytope.
\end{thm}
The Upper Bound theorem was conjectured by Motzkin in 1957
and proved by P.~McMullen in 1970.

Note that, since $C^n(m)$ is neighborly,
$$
  f_i\bigl(C^n(m)\bigr)=\bin m{i+1}\quad\text{ for }i=0,\ldots,\sbr n2-1.
$$
Due to the Dehn--Sommerville equations this determines the full $f$-vector of
$C^n(m)$.

\begin{exe}
Calculate the $f$-vector of a neighbourly simplicial $n$-polytope.
\end{exe}

\begin{exe}
\label{ubth}
Prove that the UBT inequalities are equivalent to the following inequalities
for the $h$-vector:
$$
  h_i(Q)\le\bin{m-n+i-1}i,\qquad i=0,\ldots,\sbr n2.
$$
(here $m=f_0$ is the number of vertices).
This was the key observation in McMullen's proof of the UBT and its subsequent
generalisation to arbitrary simplicial spheres due to Stanley. We will return
to this argument later.
\end{exe}

A simplicial $n$-polytope $Q$ is called \emph{stacked} if there is a
sequence $Q_0$, $Q_1,\dots$, $Q_k=Q$ of $n$-polytopes such that $Q_0$ is
an $n$-simplex and $Q_{i+1}$ is obtained from $Q_i$ by adding a pyramid
over some facet of $Q_i$. In the combinatorial language, stacked
polytopes are those obtained from
a simplex by applying several subsequent stellar subdivisions of facets.

\begin{thm}[LBT for simplicial polytopes]
\label{lbt}
  For any simplicial $n$-poly\-to\-pe $Q$ {\rm($n\ge3$)}
  with $m=f_0$ vertices the following lower bounds hold
  \begin{align*}
    f_i(Q)&\ge \bin ni f_0-\bin{n+1}{i+1}i \quad\text{ for }
    i=1,\ldots,n-2;\\
    f_{n-1}(Q)&\ge(n-1)f_0-(n+1)(n-2).
  \end{align*}
  The equality is achieved if and only if $Q$ is a stacked polytope.
\end{thm}

Note that the first inequality above is equivalent to the inequality
$h_1\le h_2$ from Theorem~\ref{gth}.

The boundary of a simplicial $n$-polytope is a simplicial subdivision
of $(n-1)$-sphere. However, not every
triangulation of sphere is combinatorially equivalent to a boundary
of polytope.

The Dehn--Sommerville equations, Upper Bound and Lower Bound theorems
were proved in different times for arbitrary triangulations of spheres.
However, the generalisation of Theorem~\ref{gth} to arbitrary simplicial
spheres remains the main combinatorial conjecture concerning the number of
faces.

\begin{prob}[$g$-conjecture]
Is it true that Theorem~\ref{gth} holds for arbitrary triangulation
of $(n-1)$-sphere?
\end{prob}

\section{Stanley--Reisner rings: combinatorics and homological algebra.}
The notion of Stanley--Reisner face ring $\k[K]$ of simplicial
complex $K$ is central in the algebraic study of triangulations.
It allows to translate the combinatorics into commutative
homological algebra.  We review its main properties, emphasising
functoriality. Then we introduce the bigraded $\Tor$-algebra
$\Tor_{\k[v_1,\dots,v_m]}(\k[K],\k)$ through a finite free
resolution of $\Z[K]$ as a module over the polynomial ring. The
corresponding bigraded Betti numbers are important combinatorial
invariants of $K$.

\subsection{Definition of $\k[K]$}
Let $\k[v_1,\ldots,v_m]$ be the graded polynomial algebra
over $\k=\Z$ or a field, \ $\deg v_i=2$.  For arbitrary subset
$\s=\{i_1,\dots,i_k\}\subset[m]$ denote by $v_\sigma$ the square-free
monomial $v_{i_1}\dots v_{i_k}$.

The \emph{face ring} (or the \emph{Stanley--Reisner ring}) of a
simplicial complex $K$ on the vertex set $[m]$ is the quotient
ring
$$
  \k[K]=\k[v_1,\ldots,v_m]/\mathcal I_K,
$$
where $\mathcal I_K$ is the homogeneous ideal generated by all
monomials $v_\sigma$ such that $\s$ is not a simplex of $K$.
The ideal $\mathcal I_K$ is called the \emph{Stanley--Reisner ideal} of~$K$.

\begin{exa}
1. Let $K$ be a 2-dimensional simplicial complex shown on Figure~\ref{figfr}.
Then
$$
  \mathcal I_K=(v_1v_5,v_3v_4,v_1v_2v_3,v_2v_4v_5).
$$
\begin{figure}[h]
\begin{center}
\begin{picture}(60,35)
\put(5,5){\line(1,0){50}}
\put(5,5){\line(1,1){25}}
\put(5,5){\line(2,1){20}}
\put(30,30){\line(1,-1){25}}
\put(25,15){\line(1,0){10}}
\put(25,15){\line(1,3){5}}
\put(35,15){\line(-1,3){5}}
\put(35,15){\line(2,-1){20}}
\put(25,25){\line(0,-1){10}}
\put(20,20){\line(0,-1){7.5}}
\put(15,15){\line(0,-1){5}}
\put(10,10){\line(0,-1){2.5}}
\put(7.5,7.5){\line(0,-1){1.25}}
\put(12.5,12.5){\line(0,-1){3.75}}
\put(17.5,17.5){\line(0,-1){6.25}}
\put(22.5,22.5){\line(0,-1){8.75}}
\put(27.5,27.5){\line(0,-1){5}}
\put(35,25){\line(0,-1){10}}
\put(40,20){\line(0,-1){7.5}}
\put(45,15){\line(0,-1){5}}
\put(50,10){\line(0,-1){2.5}}
\put(52.5,7.5){\line(0,-1){1.25}}
\put(47.5,12.5){\line(0,-1){3.75}}
\put(42.5,17.5){\line(0,-1){6.25}}
\put(37.5,22.5){\line(0,-1){8.75}}
\put(32.5,27.5){\line(0,-1){5}}
\put(5,5){\circle*{1}}
\put(25,15){\circle*{1}}
\put(35,15){\circle*{1}}
\put(30,30){\circle*{1}}
\put(55,5){\circle*{1}}
\put(3,1){1}
\put(56,1){3}
\put(29,31){2}
\put(24,11){4}
\put(34,11){5}
\end{picture}
\caption{}
\label{figfr}
\end{center}
\end{figure}
\end{exa}

\begin{exe}
Show that every square-free monomial ideal in the polynomial ring has the form
$\mathcal I_K$ for some simplicial complex~$K$.
\end{exe}

A \emph{missing face} of $K$ is a subset $\sigma\subset[m]$ such
that $\s\notin K$, but every proper subset of $\sigma$ is a
simplex of $K$.  The Stanley--Reisner ideal $\mathcal I_K$ has
basis of monomials $v_\s$ corresponding to missing faces $\s$
of~$K$.  $K$ is called a \emph{flag complex} if any set of
vertices which are pairwise connected spans a simplex of~$K$.
Alternatively, $K$ is a flag complex if and only if every its
missing face has two vertices.  The Stanley--Reisner ring $\k[K]$
is a \emph{quadratic algebra} (i.e.  the ideal $\mathcal I_K$ is
generated by quadratic monomials) if and only if $K$ is a flag
complex.

Let $K_1$ and $K_2$ be two simplicial complexes on the vertex sets $[m_1]$
and $[m_2]$ respectively. A set map $\phi\colon[m_1]\to [m_2]$ is called a
\emph{simplicial map} between $K_1$ and $K_2$ if $\phi(\sigma)\in K_2$ for
any $\sigma\in K_1$.

\begin{prop}\label{frmap}
Let $\phi\colon K_1\to K_2$ be a simplicial map. Define the map
$\phi^*\colon \k[w_1,\ldots,w_{m_2}]\to\k[v_1,\ldots,v_{m_1}]$ by
$$
  \phi^*(w_j):=\sum_{i\in\phi^{-1}(j)}v_i.
$$
Then $\phi^*$ descends to a homomorphism $\k[K_2]\to\k[K_1]$
(which we will also denote by $\phi^*$).
\end{prop}
\begin{proof}
We have to check that $\phi^*(\mathcal I_{K_2})\subset\mathcal I_{K_1}$.
Suppose
$\tau=\{j_1,\ldots,j_s\}\subset[m_2]$ is not a simplex of~$K_2$. Then
\begin{equation}
\label{fstsum}
  \phi^*(w_{j_1}\cdots w_{j_s})=\sum_{i_1\in \phi^{-1}(j_1),\ldots,
  i_s\in \phi^{-1}(j_s)}v_{i_1}\cdots v_{i_s}.
\end{equation}
We claim that $\s=\{i_1,\ldots,i_s\}$ is not a simplex of $K_1$ for any
monomial $v_{i_1}\cdots v_{i_s}$ in the right hand side of the above identity.
Indeed, otherwise we would have $\phi(\s)=\tau\in K_2$ by the definition of
simplicial map, which is impossible. Hence, the right hand side
of~\eqref{fstsum} is in~$\mathcal I_{K_1}$.
\end{proof}

Assuming that $\k$ is a field, for any graded $\k$-module $M=M^0\oplus
M^1\oplus\dots$ define its \emph{Poincar\'e series} by
$$
  F(M;t)=\sum_{i=0}^\infty(\dim_\k M^i)t^i.
$$

\begin{lem}[Stanley]
\label{psfr}
  The Poincar\'e series of $\k[K^{n-1}]$ can be calculated as
  $$
    F\bigl(\k[K^{n-1}];t\bigr)=
    \sum_{i=-1}^{n-1}\frac{f_it^{2(i+1)}}{(1-t^2)^{i+1}}=
    \frac{h_0+h_1t^2+\dots+h_nt^{2n}}{(1-t^2)^n},
  $$
  where $(f_0,\ldots,f_{n-1})$ is the $f$-vector and $(h_0,\ldots,h_n)$ is
  the $h$-vector of~$K$.
\end{lem}
\begin{proof}
A monomial in $\k[K]$ has the form $v_{i_1}^{\alpha_1}\cdots
v_{i_{k+1}}^{\alpha_{k+1}}$ where $\{i_1,\ldots,i_{k+1}\}$ is a
simplex of $K^{n-1}$ and $\alpha_1,\ldots,\alpha_{k+1}$ are some
positive integers. Thus, every $k$-simplex of $K$ contributes the
summand $\frac{t^{2(k+1)}}{(1-t^2)^{k+1}}$ to the Poincar\'e
series, which proves the first identity. The second identity is an
obvious corollary of~(\ref{hvector}).
\end{proof}

\begin{exa}
1. Let $K=\D^{n-1}$. Then $f_i=\bin{n}{i+1}$ for $-1\le i\le n-1$,
\ $h_0=1$ and $h_i=0$ for $i>0$. We have
$\k[\D^{n-1}]=\k[v_1,\ldots,v_{n}]$ and
$F(\k[\D^{n-1}];t)=(1-t^2)^{-n}$, which agrees with
Lemma~\ref{psfr}.

2. Let $K$ be the boundary of an $n$-simplex. Then $h_i=1$, \
$i=0,1,\ldots,n$, and $\k[K]=\k[v_1,\ldots,v_{n+1}]/(v_1v_2\cdots
v_{n+1})$. By Lemma~\ref{psfr},
$$
  F\bigl(\k[K];t\bigr)=\frac{1+t^2+\dots+t^{2n}}{(1-t^2)^n}.
$$
\end{exa}

\subsection{Cohen--Macaulay algebras}
Suppose now that $A=\oplus_{i\ge0}A^i$ is a finitely-generated commutative
graded algebra over $\k$. We assume that $A$ is connected ($A^0=\k$) and
has only even-degree graded components, so it is commutative as either
ordinary or graded algebra. We denote by $A_+$ the positive-degree part of
$A$ and by $\mathcal H(A_+)$ the set of homogeneous elements in $A_+$. The
\emph{Krull dimension} of $A$, denoted $\dim A$, is the maximal number of
algebraically independent elements of~$A$.

A sequence $t_1,\ldots,t_n$ of algebraically independent homogeneous
elements of $A$ is called an \emph{hsop} (homogeneous system of parameters)
if $A$ is a finitely-generated $\k[t_1,\ldots,t_n]$-module (equivalently,
$A/(t_1,\ldots,t_n)$ has finite dimension as a $\k$-vector space).

\begin{lem}[N\"other normalisation lemma]
\label{noether}
Any finitely-generated graded algebra $A$ over a field $\k$ admits
an hsop. If $\k$ is of zero characteristic and $A$ is generated by
degree-two elements, then a degree-two hsop can be chosen.
\end{lem}

A degree-two hsop is called an \emph{lsop} (linear system of parameters).

A sequence $\mb t=t_1,\ldots,t_k$ of elements of $\mathcal H(A_+)$ is
called a \emph{regular sequence} if $t_{i+1}$ is not a zero divisor in
$A/(t_1,\ldots,t_i)$ for $0\le i<k$.

\begin{exe}
Prove that any regular sequence consists of algebraically independent
elements, so it generates a polynomial subring in $A$. Show further that $\mb
t$ is a regular sequence if and only if $A$ is a \emph{free}
$\k[t_1,\ldots,t_k]$-module.
\end{exe}

Algebra $A$ is called \emph{Cohen--Macaulay} if it admits a regular hsop $\mb
t$. It follows that $A$ is Cohen--Macaulay if and only if it is free
finitely generated module over its polynomial subring.  If $\k$ is a field of
zero characteristic and $A$ is generated by degree-two elements, then one can
choose $\mb t$ to be an lsop. In this case the following formula for the
Poincar\'e series of $A$ holds
$$
  F(A;t)=\frac{F\bigl( A/(t_1,\ldots,t_n);t \bigr)}{(1-t^2)^n},
$$
where $F(A/(t_1,\ldots,t_n);t)$ is a polynomial.

A simplicial complex $K$ is called \emph{Cohen--Macaulay} (over $\k$) if its
face ring $\k[K]$ is Cohen--Macaulay.

\begin{exe}
Prove that the Krull dimension of $\k[K^{n-1}]$ equals~$n$.
Compare the above formula for $F(\k[K],\k)$ with that of Lemma~\ref{psfr}.
\end{exe}

\begin{exa}
Let $K=\partial\D^2$ be the boundary of a 2-simplex. Then
$\k[K]=\k[v_1,v_2,v_3]/(v_1v_2v_3)$. The elements $v_1,v_2\in\k[K]$ are
algebraically independent, but do not form an hsop, since
$\k[K]/(v_1,v_2)\cong\k[v_3]$ is not finite-dimensional as a $\k$-space.  On
the other hand, the elements $t_1=v_1-v_3$, $t_2=v_2-v_3$ of $\k[K]$ form an
hsop, since $\k[K]/(t_1,t_2)\cong\k[t]/t^3$. It is easy to see that $\k[K]$
is a free $\k[t_1,t_2]$-module with one 0-dimensional generator~1, one
1-dimensional generator~$v_1$, and one 2-dimensional generator~$v_1^2$. Thus,
$\k[K]$ is Cohen--Macaulay and $(t_1,t_2)$ is a regular sequence.
\end{exa}

\begin{prop}[Stanley]\label{stanmv}
If $K^{n-1}$ is a Cohen--Macaulay simplicial complex, then
$\mb h(K^{n-1})=(h_0,\ldots,h_n)$ is an $M$-vector.
\end{prop}
\begin{proof}
Let $t_1,\ldots,t_n$ be a regular sequence of degree-two
elements of~$\k[K]$. Then $A=\k[K]/(t_1,\ldots,t_n)$ is a graded
algebra generated by degree-two elements, and $\dim_\k
A^{2i}=h_i$. Now the result follows from Theorem~\ref{mvect}.
\end{proof}

For arbitrary simplex $\s\in K$ define its
\emph{link} and \emph{star} as subcomplexes
\begin{align*}
  \link_K\s&=\bigl\{\tau\in K\colon \s\cup\tau\in K,\;\s\cap\tau=\varnothing
  \bigr\};\\
  \star_K\s&=\bigl\{\tau\in K\colon \s\cup\tau\in K \bigr\}.
\end{align*}
For any vertex $v\in K$ the subcomplex $\star_K v$ can be identified with
the cone over $\link_K v$. The polyhedron $|\star_K v|$ consists of all faces
of $|K|$ that contain~$v$. We omit the subscript $K$ whenever the context
allows.

The following fundamental theorem characterises Cohen--Macaulay
com\-p\-lex\-es combinatorially.

\begin{thm}[Reisner]\label{reisner}
A simplicial complex $K$ is Cohen--Macaulay over $\k$ if and only
if for any simplex $\s\in K$ (including $\s=\varnothing$) and
$i<\dim(\link\s)$ it holds that $\widetilde{H}_i(\link\s;\k)=0$.
\end{thm}

The standard $PL$ topology techniques allow to reformulate the above
theorem in purely topological terms.

\begin{prop}[Munkres]
  $K^{n-1}$ is Cohen--Macaulay over $\k$ if and only if for arbitrary point
  $x\in|K|$ it holds that
  $\widetilde{H}_i(|K|;\k)=H_i(|K|,|K|\backslash x;\k)=0$ for $i<n-1$.
\end{prop}

\begin{cor}
A triangulation of sphere is a Cohen--Macaulay complex.
\end{cor}

Now by Theorem~\ref{stanmv} one concludes that the $h$-vector of a
simplicial sphere is an $M$-vector. This argument was used by Stanley to
extend the Upper Bound Theorem (Theorem~\ref{ubt}) to simplicial spheres.

\begin{cor}[UBT for spheres, Stanley]
\label{ubtss}
  For arbitrary simplicial $(n-1)$-sphere
  $K^{n-1}$ with $m$ vertices it holds that
  $$
    f_i(K^{n-1})\le f_i\bigl(C^n(m)\bigr) \quad\text{ for }i=1,\ldots,n-1.
  $$
\end{cor}
\begin{proof}
Due to Exercise \ref{ubth}, the UBT is equivalent to the inequalities
$$
  h_i(K^{n-1})\le\bin{m-n+i-1}i,\qquad 0\le i<{\textstyle\sbr n2}.
$$
Since $\mb h(K^{n-1})$ is an $M$-vector, there exists a graded
algebra $A=A^0\oplus A^2\oplus\dots\oplus A^{2n}$ generated by
degree-two elements such that $\dim_\k A^{2i}=h_i$
(Theorem~\ref{mvect}). In particular, $\dim_\k A^2=h_1=m-n$.
Since $A$ is generated by~$A^2$, the number $h_i$ cannot exceed
the total number of monomials of degree $i$ in $(m-n)$ variables.
The latter is exactly~$\bin{m-n+i-1}i$.
\end{proof}

\subsection{Resolutions and $\Tor$-algebras}
Let $M$ be a finitely-generated graded $\k[v_1,\ldots,v_m]$-module.
A \emph{free resolution} of $M$ is an exact sequence
\begin{equation}
\label{resol}
\begin{CD}
  \ldots @>d>> R^{-i} @>d>> \ldots @>d>> R^{-1} @>d>> R^0 @>>> M @>>> 0,
\end{CD}
\end{equation}
where the $R^{-i}$ are finitely-generated free modules and the maps $d$ are
degree-preserving. The minimal number $h$ for which there exists a free
resolution~(\ref{resol}) with $R^{-i}=0$ for $i>h$ is called the
\emph{homological dimension} of $M$ and denoted $\hd M$. By the Hilbert
syzygy theorem, $\hd M\le m$.  A resolution~(\ref{resol}) can be
written as a free {\it bigraded differential $\k$-module\/} $[R,d]$, where
$R=\bigoplus R^{-i,j}$, \ $R^{-i,j}:=(R^{-i})^j$ and $d\colon R^{-i,j}\to
R^{-i+1,j}$.  The bigraded cohomology module $H[R,d]$ has $H^{-i,k}[R,d]=0$
for $i>0$ and $H^{0,k}[R,d]=M^k$. Let $[M,0]$ be the bigraded module with
$M^{-i,k}=0$ for $i>0$, \ $M^{0,k}=M^k$, and zero differential. Then the
resolution~\eqref{resol} determined a map $[R,d]\to[M,0]$ inducing an
isomorphism in cohomology.

\begin{exe}
\label{psresol}
Show that the Poincar\'e series of $M$ can be calculated from any free
resolution~(\ref{resol}) as follows.  Suppose that $R^{-i}$ has rank $q_i$
with free generators in degrees $d_{1i},\ldots,d_{q_ii}$, \ $i=1,\ldots,h$.
Then
\begin{equation}
\label{psresol1}
  F(M;t)=(1-t^2)^{-m}\sum_{i=0}^h(-1)^i(t^{d_{1i}}+\dots+t^{d_{q_ii}}).
\end{equation}
\end{exe}

\begin{exa}[Koszul resolution]
\label{koszul}
Let $M=\k$. The $\k[v_1,\ldots,v_m]$-module structure on $\k$ is defined via
the map $\k[v_1,\ldots,v_m]\to\k$ that sends each $v_i$ to~0. Let
$\L[u_1,\ldots,u_m]$ denote the exterior algebra on $m$ generators.
Turn the tensor product $R=\L[u_1,\ldots,u_m]\otimes\k[v_1,\ldots,v_m]$
(here and below we use $\otimes$ for $\otimes_\k$) into a
{\it differential bigraded algebra\/} by setting
\begin{gather}
  \bideg u_i=(-1,2),\quad\bideg v_i=(0,2),\notag\\
  \label{diff}
  du_i=v_i,\quad dv_i=0,
\end{gather}
and requiring that $d$ be a derivation of algebras. An explicit construction
of cochain homotopy shows that $H^{-i}[R,d]=0$ for $i>0$ and
$H^0[R,d]=\k$. Since $\L[u_1,\ldots,u_m]\otimes\k[v_1,\ldots,v_m]$ is a
free $\k[v_1,\ldots,v_m]$-module, it determines a free resolution of~$\k$.
This resolution is known as the \emph{Koszul resolution}. Its expanded form
is as follows:
\begin{multline*}
0\to\L^m[u_1,\ldots,u_m]\otimes\k[v_1,\ldots,v_m]\longrightarrow\cdots\\
\longrightarrow\L^1[u_1,\ldots,u_m]\otimes\k[v_1,\ldots,v_m]\longrightarrow
\k[v_1,\ldots,v_m]\longrightarrow\k\to0,
\end{multline*}
where $\L^i[u_1,\ldots,u_m]$ is the submodule of
$\L[u_1,\ldots,u_m]$ spanned by monomials of length~$i$.
\end{exa}

Let $N$ be another module; then applying the functor
$\otimes_{\k[v_1,\ldots,v_m]}N$ to resolution $[R,d]$ we get a
homomorphism of differential modules
$$
  [R\otimes_{\k[v_1,\ldots,v_m]}N,d]\to[M\otimes_{\k[v_1,\ldots,v_m]}N,0],
$$
which in general does not induce an isomorphism in cohomology. The
$(-i)$th cohomology module of the cochain complex
$$
\begin{CD}
  \ldots @>>> R^{-i}\otimes_{\k[v_1,\ldots,v_m]}N @>>> \ldots @>>>
  R^0\otimes_{\k[v_1,\ldots,v_m]}N @>>> 0
\end{CD}
$$
is denoted $\Tor^{-i}_{\k[v_1,\ldots,v_m]}(M,N)$. Thus,
$$
  \Tor^{-i}_{\k[v_1,\ldots,v_m]}(M,N):=
  \frac{\Ker\bigl[d\colon R^{-i}\otimes_{\k[v_1,\ldots,v_m]}
  N\to R^{-i+1}\otimes_{\k[v_1,\ldots,v_m]} N\bigr]}
  {d(R^{-i-1}\otimes_{\k[v_1,\ldots,v_m]}N)}.
$$
Since all the $R^{-i}$ and $N$ are graded modules, we actually have
the \emph{bigraded} $\k$-module
$$
  \Tor_{\k[v_1,\ldots,v_m]}(M,N)=
  \bigoplus_{i,j}\Tor^{-i,j}_{\k[v_1,\ldots,v_m]}(M,N).
$$

The following properties of $\Tor^{-i}_{\k[v_1,\ldots,v_m]}(M,N)$ are
well known.
\begin{prop}
\label{torprop}
{\rm (a)} The module $\Tor^{-i}_{\k[v_1,\ldots,v_m]}(M,N)$ does not depend,
up to isomorphism, on a choice of resolution~{\rm(\ref{resol})}.

{\rm (b)} Both $\Tor^{-i}_{\k[v_1,\ldots,v_m]}(\:\cdot\:,N)$ and
$\Tor^{-i}_{\k[v_1,\ldots,v_m]}(M,\:\cdot\:)$ are covariant functors.

{\rm (c)} $\Tor^0_{\k[v_1,\ldots,v_m]}(M,N)\cong
M\otimes_{\k[v_1,\ldots,v_m]}N$.

{\rm (d)} $\Tor^{-i}_{\k[v_1,\ldots,v_m]}(M,N)\cong
\Tor^{-i}_{\k[v_1,\ldots,v_m]}(N,M)$.
\end{prop}

Now put $M=\k[K]$ and $N=\k$. Since $\deg v_i=2$, we have
$$
  \Tor_{\k[v_1,\ldots,v_m]}\bigl(\k[K],\k\bigr)=
  \bigoplus_{i,j=0}^m\Tor^{-i,2j}_{\k[v_1,\ldots,v_m]}\bigl(\k[K],\k\bigr)
$$
(i.e. $\Tor_{\k[v_1,\ldots,v_m]}(\k[K],\k)$ is non-zero only in
even second degrees). Define the \emph{bigraded Betti numbers} of
$\k[K]$ by
\begin{equation}
\label{bbnfr}
  \b^{-i,2j}\bigl(\k[K]\bigr):=
  \dim_\k\Tor^{-i,2j}_{\k[v_1,\ldots,v_m]}\bigl(\k[K],\k\bigr),\qquad
  0\le i,j\le m.
\end{equation}
We also set
$$
  \b^{-i}(\k[K])=\dim_\k\Tor^{-i}_{\k[v_1,\ldots,v_m]}(\k[K],\k)=
  \sum_j\b^{-i,2j}(\k[K]).
$$

\begin{exa}\label{ressqu}
Let $K$ be the boundary of a square. Then
$$
  \k[K]\cong\k[v_1,\ldots,v_4]/(v_1v_3,v_2v_4).
$$
Let us construct a resolution of $\k[K]$ and calculate the
corresponding bigraded Betti numbers. The module $R^0$ has one
generator 1 (of degree 0), and the map $R^0\to\k[K]$ is the
quotient projection. Its kernel is the ideal $\mathcal I_{K^1}$,
generated by two monomials $v_1v_3$ and $v_2v_4$. Take $R^{-1}$ to
be the free module on two 4-dimensional generators, denoted
$v_{13}$ and $v_{24}$, and define $d\colon R^{-1}\to R^0$ by
sending $v_{13}$ to $v_1v_3$ and $v_{24}$ to $v_2v_4$. Its kernel
is generated by one element $v_2v_4v_{13}-v_1v_3v_{24}$. Hence,
$R^{-2}$ has one generator of degree~8, say~$a$, and the map
$d\colon R^{-2}\to R^{-1}$ is injective and sends $a$ to
$v_2v_4v_{13}-v_1v_3v_{24}$. Thus, we have resolution
$$
\begin{CD}
  0 @>>> R^{-2} @>>> R^{-1} @>>> R^0 @>>> M @>>> 0,
\end{CD}
$$
where $\mathop{\rm rank}R^{0}=\b^{0,0}(\k[K])=1$,\ $\mathop{\rm
rank}R^{-1}=\b^{-1,4}=2$,\ $\mathop{\rm rank}R^{-2}=\b^{-2,8}=1$.
\end{exa}

The Betti numbers $\beta^{-i,2j}(\k[K])$ are important combinatorial
invariants of simplicial complex~$K$. The following result expresses
them in terms of homology groups of subcomplexes of~$K$.

\begin{thm}[Hochster]
\label{hoch}
We have
$$
  \beta^{-i,2j}\bigl(\k[K]\bigr)=\sum_{\omega\subset[m]\colon|\omega|=j}
  \dim_\k\widetilde{H}^{j-i-1}(K_\omega;\k),
$$
where $K_\omega$ is the full subcomplex of $K$ spanned
by~$\omega$. We assume $\widetilde{H}^{-1}(\varnothing)=\k$ above.
\end{thm}

The original Hochster's proof of this theorem uses rather complicated
combinatorial and commutative algebra techniques. Later we give
another topological proof.

\begin{exe}
Calculate the bigraded Betti numbers from Example~\ref{ressqu}
using Hochster theorem.
\end{exe}

Now let us consider the differential bigraded algebra
$[\L[u_1,\ldots,u_m]\otimes\k[K],d]$ with $d$ defined as
in~{\rm(\ref{diff})}.

\begin{lem}
\label{koscom} $\Tor_{\k[v_1,\ldots,v_m]}(\k[K],\k)$ is an algebra
in a canonical way, and there is an isomorphism of algebras
$$
  \Tor_{\k[v_1,\ldots,v_m]}(\k[K],\k)\cong
  H\bigl[\L[u_1,\ldots,u_m]\otimes\k[K],d\bigr].
$$
\end{lem}
\begin{proof}
  Using the Koszul resolution
  $[\L[u_1,\ldots,u_m]\otimes\k[v_1,\ldots,v_m],d]$ in the definition of
  $\Tor_{\k[v_1,\ldots,v_m]}(\k,\k[K])$, we calculate
  \begin{align*}
    \Tor_{\k[v_1,\ldots,v_m]}(\k[K],\k)
    &\cong\Tor_{\k[v_1,\ldots,v_m]}(\k,\k[K])\\
    &=H\bigl[ \L[u_1,\ldots,u_m]\otimes\k[v_1,\ldots,v_m]
    \otimes_{\k[v_1,\ldots,v_m]}\k[K] \bigr]\\
    &\cong H\bigl[ \L[u_1,\ldots,u_m]\otimes\k[K]\bigr].
  \end{align*}
\end{proof}

The bigraded algebra $\Tor_{\k[v_1,\ldots,v_m]}(\k[K],\k)$ is
called the \emph{$\Tor$-algebra} of simplicial complex~$K$.

\begin{lem}
\label{tamap}
A simplicial map $\phi\colon K_1\to K_2$ between two simplicial complexes on the
vertex sets $[m_1]$ and $[m_2]$ respectively induces a homomorphism
\begin{equation}\label{toralgmap}
  \phi_t^*\colon \Tor_{\k[w_1,\ldots,w_{m_2}]}\bigl(\k(K_2),\k\bigr)\to
  \Tor_{\k[v_1,\ldots,v_{m_1}]}\bigl(\k(K_1),\k\bigr)
\end{equation}
of the corresponding $\Tor$-algebras.
\end{lem}
\begin{proof}
This follows directly from Propositions~\ref{frmap}
and~\ref{torprop}~(b).
\end{proof}

\subsection{Gorenstein* complexes}
Cohen--Macaulay complexes may be characterised by means of Betti
numbers as follows.

\begin{exe}\label{cmbet}
$K^{n-1}$ is a Cohen--Macaulay complex if and only if $\b^{-i}(\k[K])=0$ for
$i>m-n$. In this case $\b^{-(m-n)}(\k[K])\ne0$.
\end{exe}

A Cohen--Macaulay complex $K^{n-1}$ on the set $[m]$ is called
\emph{Gorenstein} if $\b^{-(m-n)}(\k[K])=1$, that is,
$\Tor^{-(m-n)}_{\k[v_1,\ldots,v_m]}(\k[K],\k)\cong\k$. If additionaly $K$ is
not a cone over another simplicial complex then it is called \emph{Gorenstein*}.
In the latter case we have
$$
  \b^{-(m-n)}(\k[K])=\b^{-(m-n),2m}(\k[K])=1.
$$

The following result is similar to Reisner's theorem~\ref{reisner} and
characterises Gorenstein* simplicial complexes.

\begin{thm}[Stanley]
\label{gorencom}
  A simplicial complex $K$ is Gorenstein* over $\k$ if and only if
  for any simplex $\s\in K$ (including $\s=\varnothing$) the subcomplex
  $\link\s$ has the homology of a sphere of dimension $\dim(\link\s)$.
\end{thm}

In particular, simplicial spheres and simplicial homology spheres
(triangulated manifolds with the homology of a sphere) are Gorenstein*
complexes. However, the Gorenstein* property does not guarantee a complex to
be a triangulated manifold.  Gorenstein complexes are known to topologists as
``generalised homology spheres".

A graded commutative finite-dimensional connected $\k$-algebra
$H=\oplus_{i=0}^dH^i$ is called a \emph{Poincar\'e algebra} if $\k$-linear
maps $H^i\to\Hom_\k(H^{d-i},H^d)$, \ $a\mapsto\phi_a$, \ $\phi_a(b)=ab$ are
isomorphisms for all $i=0,\ldots,d$. The following is a corollary of
Avramov--Golod theorem~\cite[Thm. 3.4.5]{br-he98}.

\begin{thm}
A simplicial complex is Gorenstein if and only if its $\Tor$-algebra is
Poincar\'e.
\end{thm}

It is easy to see that the Poincar\'e algebra structure respects the
bigrading, whence the next result follows.

\begin{cor}
\label{tordual}
Suppose $K^{n-1}$ is a Gorenstein* complex on $[m]$. Then the following
identities hold for the Poincar\'e series of $\Tor$-algebra and for
the bigraded Betti numbers:
\begin{align*}
  &F\bigl(\Tor^{-i}_{\k[v_1,\ldots,v_m]}(\k[K],\k);\;t\bigr)=
  t^{2m}F\bigl(\Tor^{-(m-n)+i}_{\k[v_1,\ldots,v_m]}
  (\k[K],\k);\;{\textstyle\frac 1t}\bigr);\\[1mm]
  &\beta^{-i,2j}\bigl(\k[K]\bigr)=\beta^{-(m-n)+i,2(m-j)}\bigl(\k[K]\bigr),
  \quad i=0,\ldots,m-n,\;j=0,\ldots,m.
\end{align*}
\end{cor}

\begin{exe}
Deduce from the above that if $K^{n-1}$ is Gorenstein* then
$$
  F\bigl(\k[K],t\bigr)=(-1)^nF\bigl(\k[K],{\textstyle\frac1t}\bigr).
$$
Deduce further that the Dehn--Sommerville relations $h_i=h_{n-i}$, \ $0\le i\le n$,
hold for arbitrary Gorenstein* complex $K^{n-1}$.
\end{exe}

We see that the class of Gorenstein* complexes is in a sense
``the best possible algebraic approximation" to triangulated spheres.
As it was conjectured by Stanley, the $g$-theorem may continue to hold for
Gorenstein* complexes.

Later we deduce the generalised Dehn--Sommerville equations for
triangulated manifolds as a consequence of the bigraded
Poincar\'e duality for moment-angle complexes. In particular, this gives
the following short form of the equations in terms of the $h$-vector:
$$
  h_{n-i}-h_i=(-1)^i\bigl(\chi(K^{n-1})-\chi(S^{n-1})\bigr)\bin ni,
  \quad i=0,1,\ldots,n.
$$
Here $\chi(K^{n-1})=f_0-f_1+\ldots+(-1)^{n-1}f_{n-1}=1+(-1)^{n-1}h_n$ is the
Euler characteristic of $K^{n-1}$ and $\chi(S^{n-1})=1+(-1)^{n-1}$ is
that of a sphere. Note that the above equations reduce to the classical
$h_{n-i}=h_i$ in the case when $K$ is a simplicial sphere or has odd
dimension. Note also that a triangulated manifold is not Gorenstein*, or even
Cohen--Macaulay, in general.

\section[Manifolds and complexes]{Davis--Januszkiewicz spaces and moment-angle manifolds
and complexes.} Here we study different "topological models" for
the algebraic objects introduced before. The Davis--Januszkiewicz
space $\djs(K)$ has the cohomology isomorphic to the
Stanley--Reisner ring $\k[K]$. A Borel construction type model of
such space is due to Davis and Januszkiewicz~\cite{da-ja91},
whence the name comes. It also appeared in different disguises in
the work of Hattori--Masuda on torus manifolds. In our approach we
use another model for $\djs(K)$, defined through a simple colimit
(or "nested union") of nice building blocks. It is homotopy
equivalent to the Davis--Januszkiewicz model. The building block
space depends on the coefficients $\k$ (e.g., it is $\C P^\infty$
if $\k=\Z$, and $\R P^\infty$ if $\k=\Z/2$). The term
"moment-angle complex" refers to a special bigraded cellular
decomposition of the universal space $\zk$ acted on by the torus. This
space was also introduced by Davis and Januszkiewicz.
It was defined as the fibre
of a bundle with total space $\djs(K)$ and base $BT^m$, but it
has also many other interesting interpretations. The moment-angle
complex $\zk$ is a manifold provided that $K$ is a triangulation
of a sphere. At the end we discuss how the spaces $\djs(K)$ and
$\zk$ are related to toric varieties and fans,
Davis--Januszkiewicz (quasi)toric manifolds and characteristic
functions, and Hattori--Masuda torus manifolds and multifans.

\subsection{Definitions and main properties}
Here we assume $\k=\Z$, unless otherwise specified.

The classifying space for the circle $S^1$ can be identified with the
infinite-dimensional projective space $\C P^\infty$. It has a canonical cell
decomposition with one cell $B^{2k}$ in each even dimension. The classifying
space $BT^m$ of the $m$-torus is thus the product
of $m$ copies of $\C P^\infty$. The cohomology of $BT^m$
is the polynomial ring $\Z[v_1,\ldots,v_m]$, \ $\deg v_i=2$. The total
space $ET^m$ of the universal principal $T^m$-bundle over $BT^m$ can be identified
with the product of $m$ infinite-dimensional spheres.

For arbitrary subset $\omega\subset[m]$ define the subproduct
$$
  BT^\omega:=\bigl\{(x_1,\dots,x_m)\in BT^m\colon x_i=*\text{ if
  }i\notin\omega\bigr\}.
$$

For arbitrary simplicial complex $K$ on $[m]$ define the \emph{Davis-Januszkiewicz space}
as the following cellular subcomplex:
$$
  \djs(K):=\bigcup_{\s\in K}BT^\s \subset BT^m.
$$

\begin{prop}
\label{homsrs} The cohomology of $\djs(K)$ is isomorphic to the
Stanley--Reisner ring~$\Z[K]$. Moreover, the inclusion of cellular
complexes $i\colon \djs(K)\hookrightarrow BT^m$ induces the quotient
epimorphism
$$
  i^*\colon \Z[v_1,\ldots,v_m]\to\Z[K]=\Z[v_1,\ldots,v_m]/\mathcal I_K
$$
in the cohomology.
\end{prop}
\begin{proof}
A monomial $v_{i_1}^{k_1}\dots v_{i_p}^{k_p}$ represents the cellular cochain
$(B^{2k_1}_{i_1}\dots B^{2k_p}_{i_p})^*$ in $C^*(BT^m)$. Under the cochain homomorphism
induced by the inclusion $\djs(K)\subset BT^m$ it maps identically if $\{i_1,\dots,i_p\}\in K$
and to zero otherwise, whence the statement follows.
\end{proof}

It is convenient to realise the torus $T^m$ as a subspace in $\C^m$:
$$
  T^m=\bigl\{ (z_1,\ldots,z_m)\in\C^m\colon |z_i|=1,\quad i=1,\ldots,m
  \bigr\}.
$$
It is contained in the \emph{unit polydisc}
$$
  (D^2)^m:=\bigl\{ (z_1,\ldots,z_m)\in\C^m\colon |z_i|\le1,\quad i=1,\ldots,m
  \bigr\}.
$$
For arbitrary subset $\omega\subset[m]$ define
$$
  B_\omega:=\bigl\{ (z_1,\ldots,z_m)\in(D^2)^m\colon
  |z_i|=1\text{ if }i\notin\omega\bigr\}.
$$
Obviously, $B_\omega$ is homeomorphic to $(D^2)^{|\omega|}\times T^{m-|\omega|}$.
Given a simplicial complex $K$ on $[m]$, define the \emph{moment-angle complex} $\zk$ by
\begin{equation}\label{bdeco}
  \zk:=\bigcup_{\s\in K}B_\s \subset (D^2)^m.
\end{equation}

\begin{rem}
The above constructions of $\djs(K)$ and $\zk$ are examples of
\emph{colimits} of diagrams of topological spaces over the
\emph{face category} of $K$ (objects are simplices and morphisms
are inclusions). In the first case the diagram assigns the space
$BT^\s$ to a simplex $\s$; its colimit is $\djs(K)$. The second
diagram assigns $B_\s$ to $\s$; its colimit is $\zk$.
\end{rem}

The torus $T^m$ acts on $(D^2)^m$ coordinatewise, and each subspace $B_\omega$ is
invariant under this action. Therefore, the space $\zk$ itself is acted on by
the torus. The quotient $(D^2)^m/T^m$ can be identified with the \emph{unit $m$-cube}:
$$
  I^m:=\bigl\{ (y_1,\ldots,y_m)\in\R^m\colon 0\le y_i\le1,\quad i=1,\ldots,m
  \bigr\}.
$$
The quotient $B_\omega/T^m$ is then the following $|\omega|$-face of $I^m$:
$$
  C_\omega:=\bigl\{ (y_1,\ldots,y_m)\in I^m\colon y_i=1\text{ if
  }i\notin\omega\bigr\}.
$$
Thus, the whole quotient $\zk/T^m$ is identified
with a certain cubical subcomplex in $I^m$, which we denote $\cc(K)$.

\begin{lem}
The cubical complex $\cc(K)$ is $PL$-homeomorphic to $|\cone K|$.
\end{lem}
\begin{proof}
Let $K'$ denote the barycentric subdivision of $K$ (the vertices of $K'$ correspond to
non-empty simplices $\s$ of $K$). We define a $PL$ embedding $i_c\colon \cone K'\hookrightarrow I^m$
by sending each vertex $\s$ to the vertex $(\e_1,\dots,\e_m)\in I^m$, where $\e_i=0$ if $i\in\s$
and $\e_i=1$ otherwise, sending the cone vertex to $(1,\dots,1)\in I^m$, and extending
linearly on the simplices of $\cone K'$. The barycentric subdivision of a
face $\s\in K$ is a
subcomplex in $K'$, which we denote $K'|_\s$. Under the map $i_c$ the subcomplex
$\cone K'|_\s$ maps onto the face $C_\s\subset I^m$. Thus, the whole complex $\cone K'$
maps homeomorphically onto $\cc(K)$, whence the proof follows.
\end{proof}

\begin{exa}
If $K=\D^{m-1}$ is the whole simplex on $[m]$, then $\cc(K)$ is the whole cube $I^m$,
and the above constructed $PL$-homeomorphism between $\cone(\D^{m-1})'$ to $I^m$ defines
the \emph{canonical triangulation} of $I^m$.
\end{exa}

\begin{lem}
Suppose that $K$ is a triangulation of a sphere: $|K|\cong S^{n-1}$. Then $\zk$
is an $(m+n)$-dimensional manifold.
\end{lem}
\begin{proof}
The space $|\cone K'|$ has a canonical \emph{face structure} whose
facets (co\-di\-men\-si\-on-\-one faces) are
$F_i:=\star_{K'}\{i\}$, \ $i=1,\dots,m$, and $i$-faces are
non-empty intersections of $i$-tuples of facets. Since $K$ is a
triangulation of a sphere, $|\cone K'|$ is an $n$-ball. Every
point in $|\cone K'|$ has a neighbourhood homeomorphic to an open
subset in $I^n$ (or $\R^n_+$) with the homeomorphism preserving
the dimension of faces. By the definition, this displays $|\cone
K'|$ as a \emph{manifold with corners}. Having identified $|\cone
K'|$ with $\cc(K)$ and further $\cc(K)$ with $\zk/T^m$, we see
that every point in $\zk$ lies in a neighbourhood homeomorphic to
an open subset in $(D^2)^n\times T^{m-n}$ and thus in $\R^{m+n}$.
\end{proof}

Suppose $P^n$ is a \emph{simple} $n$-polytope, that is,
every vertex of $P$ is contained in exactly $n$ facets. Define $K_P$ to be
the boundary complex of the dual simplicial polytope. Then the face structure
in $\cone K'_P$ is the same as that of $P^n$.

\begin{exe}
Show that the isotropy subgroup of a point $x\in\cc(K)$ under the $T^m$-action on $\zk$
is the coordinate subtorus
$$
  T(x):=\bigl\{ (z_1,\ldots,z_m)\in T^m\colon z_i=1\text{ if }x\notin F_i
  \bigr\}.
$$
In particular, the action is free over the interior (that is,
around the cone point) of $\cc(K)\cong|\cone K'|$.
\end{exe}

It follows that the moment-angle complex can be identified with the
quotient
$$
  \zk=\bigl(T^m\times|\cone K'|\bigr)/{\sim},
$$
where $(t_1,x)\sim(t_2,y)$ if and only if $x=y$ and $t_1t_2^{-1}\in T(x)$.
In the case $K=K_P$ we may write $(T^m\times P^n)/{\sim}$
instead. The latter $T^m$-manifold is the one introduced by Davis and
Januszkiewicz~\cite{da-ja91},
which thereby coincides with our moment-angle complex.

For arbitrary $T^m$-space $X$ define the \emph{Borel
construction} (also known as the \emph{homotopy quotient} or \emph{associated
bundle}) as the identification space
$$
  ET^m\times_{T^m} X:=ET^m\times X/{\sim},
$$
where $(e,x)\sim(eg,g^{-1}x)$ for any $e\in ET^m$, \ $x\in X$, $g\in T^m$.

The projection $(e,x)\to e$ identifies $ET^m\times_{T^m} X$ with the total space
of a bundle $ET^m\times_{T^m} X\to BT^m$ with fibre $X$ and structure
group~$T^m$.

In the sequel we denote the Borel construction $ET^m\times_{T^m} X$
by~$B_TX$.

In particular, for any simplicial complex $K$ on $m$ vertices we
have the Borel construction $B_T\zk$ and the bundle $p\colon B_T\zk\to
BT^m$ with fibre~$\zk$. Davis and Januszkiewicz showed that the
cohomology of $B_T\zk$ is isomorphic to the Stanley--Reisner
ring~$\Z[K]$, and $p^*\colon \Z[v_1,\ldots,v_m]\to\Z[K]$ is the quotient
projection. Thus, the space $B_T\zk$ provides another topological
model for the Stanley--Reisner ring. However, the following result
shows that the two models $\djs(K)$ and $B_T\zk$ are homotopy
equivalent.

\begin{thm}
\label{homeq1}
  There is a deformation retraction $B_T\zk\to\djs(K)$ such that the diagram
  $$
  \begin{CD}
    B_T\zk @>p>> BT^m\\
    @VVV @|\\
    \djs(K) @>i>> BT^m
  \end{CD}
  $$
  is commutative.
\end{thm}
\begin{proof}
We have $\zk=\bigcup_{\s\in K}B_\s$, and each $B_\s$ is $T^m$-invariant.
Hence, there is the corresponding decomposition of the Borel construction:
$$
  B_T\zk=ET^m\times_{T^m}\zk=\bigcup_{\s\in K}ET^m\times_{T^m}B_\s.
$$
Suppose $|\s|=s$. Then $B_\s\cong(D^2)^s\times T^{m-s}$, so we have
$ET^m\times_{T^m}B_\s\cong(ET^s\times_{T^s}(D^2)^s)\times ET^{m-s}$.
The space
$ET^s\times_{T^s}(D^2)^s$ is the total space of a $(D^2)^s$-bundle
over $BT^s$, and $ET^{m-s}$ is contractible. It follows that there is a deformation retraction
$ET^m\times_{T^m}B_\s\to BT^\s$. These homotopy equivalences
corresponding to different simplices fit together to yield a
required homotopy equivalence between $p\colon B_T\zk\to BT^m$ and
$i\colon \djs(K)\hookrightarrow BT^m$.
\end{proof}

Below we denote by $\djs(K)$ either of the two homotopy equivalent spaces.

\begin{cor}
\label{homfib}
  The moment-angle complex $\zk$ is the homotopy fibre of the embedding
  $i\colon \djs(K)\hookrightarrow BT^m$.
\end{cor}

\begin{cor}\label{eqczk}
The $T^m$-equivariant cohomology of $\zk$ is isomorphic to the
Stanley--Reisner ring of~$K$:
$$
  H^*_{T^m}(\zk)\cong\Z[K].
$$
\end{cor}

The following information about the homotopy groups of $\zk$ can be retrieved
from the above constructions.

\begin{prop}
\label{homgr}
  {\rm(a)} $\zk$ is {\rm2}-connected, and $\pi_i(\zk)=\pi_i(\djs(K))$ for
  $i\ge 3$.

  {\rm(b)}  If $K$ is $q$-neighbourly,
  then $\pi_i(\zk)=0$ for $i<2q+1$. Moreover,
  $\pi_{2q+1}(\zk)$ is a free Abelian group generated by
  the $(q+1)$-element missing faces of $K$.
\end{prop}
\begin{proof}
Note that $BT^m=K(\Z^m,2)$ and the 3-skeleton of $\djs(K)$
coincides with that of~$BT^m$. If $K$ is $q$-neighbourly, then the
$(2q+1)$-skeleton of $\djs(K)$ coincides with that of~$BT^m$. Now,
both statements follow easily from the exact homotopy sequence of
the map $i\colon \djs(K)\to BT^m$.
\end{proof}

\subsection{Cell decompositions}
Here we construct a canonical cell decomposition of the moment-angle complex.

\begin{figure}[h]
  \begin{picture}(120,35)
  \put(60,20){\oval(30,30)}
  \put(75,20){\circle*{1.5}}
  \put(76,21){1}
  \put(42,21){$T$}
  \put(52,25){$D$}
  \end{picture}
  \caption{}
  \label{cellfig}
\end{figure}
Let us decompose $D^2$ into the union of one 2-cell~$D$, one
1-cell $T$ and one 0-cells~1, as shown on Figure~\ref{cellfig}. It
defines a cellular complex structure on the polydisc $(D^2)^m$ with $3^m$ cells.
Each cell of this complex is a product of cells of 3 different
types: $D_i$, $T_i$ and $1_i$, \ $i=1,\ldots,m$. We encode
the cells using the language of ``sign vectors".
Each cell of $(D^2)^m$ will be represented by a sign vector $\mathcal T\in\{D,T,1\}^m$. We denote
by $\mathcal T_D$, $\mathcal T_T$ and $\mathcal T_1$ respectively
the $D$-, $T$- and $1$-component of~$\mathcal T$. Each of these
components can be seen as a subset of~$[m]$, and all three subsets are
complementary. The closure of a cell $\mathcal T$ is homeomorphic to a product
of $|\mathcal T_D|$ discs and
$|\mathcal T_T|$ circles.

\begin{lem}
\label{zkcell}
  $\zk$ is a cellular subcomplex of $(D^2)^m$.
  A cell $\mathcal T\subset(D^2)^m$
  belongs to $\zk$ if and only if $\mathcal T_D\in K$.
\end{lem}
\begin{proof}
We have $\zk=\cup_{\s\in K}B_\s$ and each $B_\s$ is the closure of
the cell $\mathcal T$ with $\mathcal T_D=\s$, \ $\mathcal T_T=[m]\setminus\s$
and $\mathcal T_1=\varnothing$.
\end{proof}

Denote by $C^*(\zk)$ the corresponding cellular cochains. It has a natural
bigrading defined by $\bideg\mathcal T=(-|\mathcal T_T|,2|\mathcal T_T|+2|\mathcal T_D|)$
(so $\bideg D_i=(0,2)$, $\bideg T_i=(-1,2)$ and $\bideg 1_i=(0,0)$). Moreover,
since the cellular differential does not change the second grading, $C^*(\zk)$
splits into the sum of its components with fixed second degree:
$$
  C^*(\zk)=\bigoplus_{j=1}^m C^{*,2j}(\zk).
$$
Correspondingly, the cohomology of $\zk$ acquires an additional
grading, and one may define the \emph{bigraded Betti numbers}
$b^{-i,2j}(\zk)$ by
$$
  b^{-i,2j}(\zk):=\dim_\k H^{-i,2j}(\zk), \quad i,j=1,\dots,m.
$$
For the ordinary Betti numbers one has
$b^k(\zk)=\sum_{2j-i=k}b^{-i,2j}(\zk)$.

\begin{lem}
\label{mamap} Let $\phi\colon K_1\to K_2$ be a simplicial map between
complexes on the sets $[m_1]$ and $[m_2]$ respectively. Then there
is an equivariant cellular map $\phi_{\mathcal Z}\colon\mathcal
Z_{K_1}\to\mathcal Z_{K_2}$ covering the induced map $|\cone
K'_1|\to|\cone K'_2|$.
\end{lem}
\begin{proof}
Define the map $\phi_D\colon (D^2)^{m_1}\to(D^2)^{m_2}$ by
$$
  \phi_D(z_1,\dots,z_{m_1})=(w_1,\dots,w_{m_2}),
$$
where
$$
  w_i=\prod_{j\in f^{-1}(i)}z_j.
$$
Then one easily checks that $\phi_D(B_\s)\subset B_{\phi(\s)}$. Since $\phi$
is simplicial, $\s\in K_1$ implies $\phi(\s)\in K_2$. Hence, the restriction
of $\phi_D$ to $\mathcal Z_{K_1}$ is the required map.
\end{proof}

The above constructed map regards the bigrading, so the bigraded Betti
numbers are functorial.

\subsection{Toric varieties, quasitoric manifolds, and torus manifolds}
Several important classes of manifolds with torus action emerge as
the quotients of moment-angle complexes by appropriate freely
acting subtori.

First we give the following characterisation of lsop's in the
Stanley--Reisner ring. Let $K^{n-1}$ be a simplicial complex and
$\t_1,\dots,\t_n$ a sequence of degree-two elements in $\k[K]$. We
may write
\begin{equation}\label{lsop}
  \t_i=\l_{i1}v_1+\dots+\l_{im}v_m,\quad i=1,\dots,n.
\end{equation}
For arbitrary simplex $\s\in K$ we have $K_\s=\D^{|\s|-1}$ and
$\k[K_\s]$ is the polynomial ring $\k[v_i\colon i\in\s]$ on $|\s|$
generators. The inclusion $K_\s\subset K$ induces the
\emph{restriction homomorphism} $r_\s$ from $\k[K]$ to the
polynomial ring, mapping $v_i$ identically if $i\in\s$ and to zero
otherwise. The degree-two part of a polynomial ring on $q$
generators may be identified with the space of linear forms
on~$\k^q$.

\begin{lem}\label{restr}
A degree-two sequence $\t_1,\dots,\t_n$ is an lsop in
$\k[K^{n-1}]$ if and only if its image under any restriction
homomorphism $r_\s$ generates $(\k^{|\s|})^*$.
\end{lem}
\begin{proof}
Suppose \eqref{lsop} is an lsop. For simplicity we denote its image under
any restriction homomorphism by the same letters. Then the restriction
induces a homomorphism of quotient rings:
$$
  \k[K]/(\t_1,\dots,\t_n)\to\k[v_i\colon i\in\s]/(\t_1,\dots,\t_n).
$$
Since \eqref{lsop} is an lsop, $\k[K]/(\t_1,\dots,\t_n)$ is a
finitely generated $\k$-module. Hence, so is
$\k[v_i\colon i\in\s]/(\t_1,\dots,\t_n)$. But the latter can be finitely
generated only if $\t_1,\dots,\t_n$ generates the degree-two part
of the polynomial ring.

The ``if" part may be proved by considering the sum of restrictions:
$$
  \k[K]\to\bigoplus_{\s\in K}\k[v_i\colon i\in\s],
$$
which is actually a monomorphism. See \cite[Thm.~5.1.16]{br-he98}
for details.
\end{proof}

In particular, if $K^{n-1}$ is pure (i.e., all maximal simplices have the same dimension),
then \eqref{lsop} is an lsop if and only if its restriction to every $(n-1)$-simplex
is a basis in the space of linear forms.

Suppose now that $K$ is Cohen--Macaulay (e.g., $K$ is a sphere triangulation). Then
every lsop is a regular sequence (however, for $\k=\Z$ or a field of finite characteristic
an lsop may fail to exist).

Now we restrict to the case $\k=\Z$ and
organise the coefficients in \eqref{lsop} into an $n\times m$-matrix $\L=(\l_{ij})$.
For arbitrary maximal simplex $\s\in K$ denote by $\L_\s$ the square submatrix (minor)
formed by the elements $\l_{ij}$ with $j\in\s$. The matrix $\Lambda$ defines a linear map
$\Z^m\to\Z^n$ and a homomorphism $T^m\to T^n$. We denote both by~$\l$ and denote
the kernel of the latter map by $T_\L$.

\begin{thm}\label{zkquo}
The following conditions are equivalent:
\begin{itemize}\label{quot}
\item[(a)] the sequence \eqref{lsop} is an lsop in $\Z[K^{n-1}]$;
\item[(b)] $\det\L_\s=\pm1$ for every maximal simplex $\s\in K$;
\item[(c)] $T_\L\cong T^{m-n}$ and $T_\L$ acts freely on $\zk$.
\end{itemize}
\end{thm}
\begin{proof}
The equivalence of (a) and (b) is a reformulation of Lemma~\ref{restr}.
Let us prove the equivalence of (b) and (c). Every isotropy subgroup
of the $T^m$-action on $\zk$ has the form
$$
  T^\s=\bigl\{ (z_1,\dots,z_m)\in T^m\colon z_i=1\text{ if }i\notin\s \bigr\}
$$
for some simplex $\s\in K$. Now, (b) is equivalent to the condition
$T_\L\cap T^\s=\{e\}$ for arbitrary maximal $\s$, whence the statement follows.
\end{proof}

We denote the quotient $\zk/T_\L$ by $M^{2n}_K(\L)$ (when the
context allows we may abbreviate it to $M^{2n}_K$ or even to
$M^{2n}$). If $K$ is a triangulated sphere, then $\zk$ is a
manifold, hence, so is $M^{2n}_K$. The $n$-torus
$T^n=T^m/T_\L$ acts on $M^{2n}_K$.
This construction generalises the following two
important classes of $T^n$-manifolds.

Let $K=K_P$ be a polytopal triangulation, dual to the boundary
complex of a simple polytope $P^n$. Then the map $\l$ determined
by the matrix $\L$ may be regarded as an assignment of an integer
vector to any facet of~$P^n$. A map $\l$ coming from a matrix
satisfying the condition of Theorem~\ref{quot}(b) was called a
\emph{characteristic map} by Davis and Januszkiewicz. We refer to
the corresponding quotient $M^{2n}_P(\L)=\mathcal Z_{K_P}/T_\L$ as
a \emph{quasitoric manifold} (a toric manifold in the
Davis--Januszkiewicz terminology).

Let us assume further that $P^n$ is realised in $\R^n$
with integer vertex coordinates, so we can write
\begin{equation}\label{ptope}
  P^n=\bigl\{\mb x\in\R^n\colon\<\mb l_i,\mb x\>\ge-a_i,
  \; i=1,\dots,m\bigr\},
\end{equation}
where $\mb l_i$ are integral vectors normal to the facets of $P^n$
(we may assume these vectors to be primitive and inward pointing), and
$a_i\in\Z$. Let $\Lambda$ be the matrix formed by the column vectors
$\mb l_i$, \ $i=1,\dots,m$. Then $\mathcal Z_{K_P}/T_\L$ is the
\emph{projective toric variety} determined by the polytope~$P^n$.
The condition of Theorem~\ref{quot}(b) is equivalent to that
the toric variety is non-singular. Thereby a non-singular projective toric
variety is a quasitoric manifold (but there are many quasitoric manifolds
which are not toric varieties).

Finally, we mention that if $K$ is an arbitrary triangulation of
sphere, then the manifold $M^{2n}_K(\L)$ is a \emph{torus manifold}
in the sense of Hattori--Masuda~\cite{ha-ma??}. The corresponding multi-fan
has $K$ as the underlying simplicial complex. This particular class of torus
manifolds has many interesting specific properties.

\section{Cohomology of moment-angle complexes}
Studying ordinary and $T^m$-equivariant topology of $\zk$ opens the
way to some combinatorial applications. While the cohomology of
$\djs(K)$ is $\Z[K]$, the cohomology of the moment-angle complex
$\zk$ is isomorphic to the Tor-algebra
$\Tor_{\k[v_1,\ldots,v_m]}(\k[K],\k)$. The argument uses some
algebraic topology techniques, such as Eilenberg-Moore spectral
sequences (shortly \emph{emss}). The intrinsic bigrading in Tor is
exactly that coming from the bigraded cell decomposition of $\zk$.
In the case when $K$ is a triangulated manifold, the bigraded relative
Poincar\'e duality for $\zk$ gives the generalised
Dehn-Sommerville equations.

\subsection{Eilenberg--Moore spectral sequence}
The Eilenberg--Moore spectral sequence can be considered as an
extension of Adams' cobar construction approach to calculating the
cohomology of loop spaces. Here we give the necessary information
on the spectral sequence; we follow L.~Smith's paper~\cite{smit67}
in this description. The following theorem provides an algebraic
setup for the emss.

\begin{thm}[Eilenberg--Moore]
\label{algemss} Let $A$ be a differential graded $\k$-algebra, and
$M$,~$N$ differential graded $A$-modules. Then there exists a
second quadrant spectral sequence $\{E_r,d_r\}$ converging to
$\Tor_A(M,N)$ whose $E_2$-term is
$$
  E_2^{-i,j}=\Tor^{-i,j}_{H[A]}\bigl(H[M],H[N]\bigr),\quad i,j\ge0,
$$
where $H[\cdot]$ denotes the cohomology algebra (or module).
\end{thm}

Topological applications of the above theorem arise in the case
when $A,M,N$ are singular cochain algebras (or their commutative
models) of certain topological spaces. The classical situation is
described by the commutative diagram
\begin{equation}
\begin{CD}
  E @>>> E_0\\
  @VVV @VVV\\
  B @>>> B_0,
\end{CD}
\label{comsq}
\end{equation}
where $E_0\to B_0$ is a Serre fibre bundle with fibre $F$ over a
simply connected base~$B_0$, and $E\to B$ is the pullback along a
continuous map $B\to B_0$. For any space~$X$, let $C^*(X;\k)$
denote the singular cochain algebra of $X$ or its commutative
model (e.g., the cellular cochain algebra if a nice cellular
product exists). Obviously, $C^*(E_0;\k)$ and $C^*(B;\k)$ are
$C^*(B_0;\k)$-modules. Under these assumptions the following
statement holds.

\begin{lem}[Eilenberg--Moore]\label{gentoralg}
$\Tor_{C^*(B_0;\k)}(C^*(E_0;\k),C^*(B;\k))$ is an algebra in a
natural way, and there is a canonical isomorphism of algebras
$$
  \Tor_{C^*(B_0;\k)}\bigl(C^*(E_0;\k),C^*(B;\k)\bigr)\to H^*(E;\k).
$$
\end{lem}

The above two results lead to the following statement.
\begin{thm}[{Eilenberg--Moore}]
There exists a spectral sequence of $\k$-\-al\-geb\-ras
$\{E_r,d_r\}$ with
\begin{enumerate}
\item[$(a)$] $E_r\Rightarrow H^*(E;\k)$;
\item[$(b)$]
$E_2^{-i,j}=\Tor^{-i,j}_{H^*(B_0;\k)}\bigl(H^*(E_0;\k),H^*(B;\k)\bigr)$.
\end{enumerate}
\end{thm}

The case when $B$ above is a point is particularly important for
applications, so we state the corresponding result separately.
\begin{cor}
  Let $E\to B$ be a fibration over a simply connected space $B$ with
  fibre~$F$. Then there exists a spectral sequence of $\k$-algebras
  $\{E_r,d_r\}$ with
  \begin{enumerate}
    \item[$(a)$] $E_r\Rightarrow H^*(F;\k)$;
    \item[$(b)$] $E_2=\Tor_{H^*(B;\k)}\bigl(H^*(E;\k),\k\bigr)$.
  \end{enumerate}
\end{cor}
We refer to the above spectral sequence as the
\emph{Eilenberg--Moore spectral sequence of fibration $E\to B$}.

Here is the first application of the emss. We use the notation of
Theorem~\ref{zkquo}.

\begin{thm}\label{cohmk}
The cohomology of the quotient $M_K(\L)=\zk/T_\L$ is given by
$$
  H^*\bigl(M_K(\L)\bigr)\cong\Z[K]/(\t_1,\dots,\t_n).
$$
\end{thm}
\begin{proof}
Denote $M:=M_K(\L)$. Since $T_\L$ acts freely on $\zk$, we have
$$
  ET^m\times_{T^m}\zk\simeq ET^n\times_{T^n}M,
$$
where $T^n=T^m/T_\L$. Hence, $H^*(ET^n\times_{T^n}M)\cong\Z[K]$ by
Corollary~\ref{eqczk}. The emss of the bundle
$ET^n\times_{T^n}M\to BT^n$ converges to $H^*(M)$ and has
$$
  E_2^{*,*}=\Tor^{*,*}_{H^*(BT^n)}\bigl(H^*(ET^n\times_{T^n}M),\Z\bigr)=
  \Tor^{*,*}_{\Z[t_1,\ldots,t_n]}\bigl(\Z[K],\Z\bigr),
$$
where the $\Z[t_1,\ldots,t_n]$-module structure in $\Z[K]$ is
defined through the homomorphism $\Z[t_1,\ldots,t_n]\to\Z[K]$
taking $t_i$ to $\theta_i$, \ $i=1,\ldots,n$, see~\eqref{lsop}.
Since $\Z[K]$ is Cohen--Macaulay and $\t_1,\dots,\t_n$ is an lsop,
$\Z[K]$ is a free $\Z[t_1,\dots,t_n]$-module, so we have
\begin{align*}
  \Tor^{*,*}_{\Z[t_1,\ldots,t_n]}\bigl(\Z[K],\Z\bigr)&=
  \Tor^{0,*}_{\Z[t_1,\ldots,t_n]}\bigl(\Z[K],\Z\bigr)\\
  &=\Z[K]\otimes_{\Z[t_1,\ldots,t_n]}\Z=\Z[K]\bigr/(\t_1,\dots,\t_n).
\end{align*}
Therefore, $E_2^{0,*}=\Z[K]/(\t_1,\dots,\t_n)$ and $E_2^{-p,*}=0$
for $p>0$. It follows that the Eilenberg--Moore spectral sequence
collapses at the $E_2$ term and $H^*(M)=\Z[K]/(\t_1,\dots,\t_n)$,
as claimed.
\end{proof}

In the case of quasitoric manifold the cohomology ring was described
in~\cite{da-ja91}, while the case of non-singular toric varieties was
known long before as the \emph{Danilov--Jurkiewicz theorem}.

\begin{cor}
The cohomology of $M_K(\L)$ vanishes in odd dimensions and
$H^{2i}(M_K(\L))$ is a free abelian group of rank $h_i$, \ $i=0,\dots,n$.
\end{cor}

In the case when $P^n$ is an integral polytope~\eqref{ptope} and
$M^{2n}_P$ a non-singular toric variety, the cohomology ring
$H^*(M_P)$ has more specific properties. Namely there is a
cohomology class $\omega\in H^2(M_P)$ for which the multiplication maps
$$
  H^{n-i}(M_P) \stackrel{\omega^i}{\longrightarrow} H^{n+i}(M_P),\qquad
  i=1,\dots,n
$$
are isomorphisms. This fact is known as the \emph{Hard Lefschetz theorem for
toric varieties} (although its proof is far beyond the scope of this review).
One can take $\omega=a_1v_1+\cdots+a_mv_m$ in the notation of~\eqref{ptope}.
The multiplication by $\omega$ determines a monomorphism
$H^{2i}(M_P)\to H^{2i+2}(M_P)$ for $i\le \sbr n2$. Consider the graded algebra
$A:=H^*(M_P)/(\omega)$. It is generated in degree two and has $\dim
A^{2i}=h_i-h_{i-1}$, \ $i=1,\dots,\sbr n2$. Then it follows from
Theorem~\ref{mvect} that
$$
  0\le h_{i+1}-h_i \le (h_i-h_{i-1})^{\langle i\rangle}
$$
for $i<\sbr n2$, which proves the necessity part of Theorem~\ref{gth}.

\subsection{Cohomology of $\zk$}
Here we calculate the cohomology algebra of $\zk$. As an immediate
corollary we obtain that the cohomology inherits a canonical
\emph{bi}grading from the spectral sequence. The corresponding
bigraded Betti numbers coincide with the algebraic Betti numbers
of~$\k[K]$ introduced before~\eqref{bbnfr}.

\begin{thm}\label{cohom1}
There are isomorphisms of algebras
$$
  H^*(\zk;\k)\cong\Tor_{\k[v_1,\dots,v_m]}\bigl(\k[K],\k\bigr)\cong
  H\bigl[\k[K]\otimes\Lambda[u_1,\dots,u_m],d\bigr],
$$
where the latter algebra is bigraded by $\bideg v_i=(0,2)$, \
$\bideg u_i=(-1,2)$ and the differential is determined by
$dv_i=0$, \ $du_i=v_i$. In particular,
$$
  \dim_\k H^q(\zk;\k)\cong
  \sum_{-i+2j=q}\beta^{-i,2j}\bigl(\k(K)\bigr).
$$
\end{thm}
\begin{proof}
Consider the emss of the commutative square
$$
\begin{CD}
  E @>>> ET^m\\
  @VVV @VVV\\
  \djs(K) @>i>> BT^m
\end{CD},
$$
where the left vertical arrow is the pullback along~$i$.
Corollary~\ref{homfib} shows that $E$ is homotopy equivalent
to~$\zk$.

By Proposition~\ref{homsrs}, the map $i\colon \djs(K)\hookrightarrow
BT^m$ induces the quotient epimorphism
$i^*\colon \k[v_1,\ldots,v_m]\to\k[K]$ of the cellular cochain algebras.
Since $ET^m$ is contractible, there is a chain equivalence
$C^{*}(ET^m;\k)\simeq\k$. More precisely, $C^*(ET^m;\k)$ can be
identified with the Koszul resolution
$\L[u_1,\ldots,u_m]\otimes\k[v_1,\ldots,v_m]$ of $\k$ (see
Example~\ref{koszul}). Therefore, we have an isomorphism
$$
 \Tor_{C^{*}(BT^m)}\bigl(C^{*}\bigl(\djs(K)\bigr),C^{*}(ET^m)\bigr)\cong
 \Tor_{\k[v_1,\ldots,v_m]}\bigl(\k[K],\k\bigr).
$$
On the other hand, Lemma~\ref{gentoralg} shows that
$\Tor_{C^{*}(BT^m)}(C^{*}(\djs(K)),C^{*}(ET^m))$ is an algebra
isomorphic to $H^{*}(\zk;\k)$, which proves the first isomorphism.
The second one follows from Lemma~\ref{koscom}.
\end{proof}

\begin{exe}
Prove that the Leray--Serre spectral sequence of the principal
$T^m$-bundle $\zk\to\djs(K)$ collapses at the $E_3$ term.
\end{exe}

\begin{exa}
Let $K=\partial\D^{m-1}$. Then
$\k[K]=\k[v_1,\ldots,v_m]/(v_1\cdots v_m)$. A direct calculation
shows that the cohomology of $\k[K]\otimes\Lambda[u_1,\ldots,u_m]$
is additively generated by the classes $1$ and $[v_1v_2\cdots
v_{m-1}u_m]$. The latter represents the fundamental cohomology
class of $\zk\cong S^{2m-1}$.
\end{exa}

\begin{exa}
Let $K$ be the boundary complex of an $m$-gon $P^2$ with $m\ge4$.
We have $\k[K]=\k[v_1,\ldots,v_m]/\mathcal I_K$, where $\mathcal
I_K$ is generated by the monomials $v_iv_j$, \ $i-j\ne0,1\mod m$.
The complex $\zk$ is a manifold of dimension~$m+2$. For $m=4$ we
have $\zk=S^3\times S^3$. Suppose $m=5$. Then the group $H^3(\zk)$
has 5 generators represented by the cocycles
$v_iu_{i+2}\in\k[K]\otimes\Lambda[u_1,\ldots,u_5]$, \
$i=1,\ldots,5$, and $H^4(\zk)$ has 5 generators represented by the
cocycles $v_ju_{j+2}u_{j+3}$, \ $j=1,\ldots,5$. The product of
$v_iu_{i+2}$ and $v_ju_{j+2}u_{j+3}$ represents a non-zero
cohomology class in $H^7(\zk)$ if and only if all the indices
$i,i+2,j,j+2,j+3$ are different. Thus, for each of the 5
cohomology classes $[v_iu_{i+2}]$ there is a unique (Poincar\'e
dual) cohomology class $[v_ju_{j+2}u_{j+3}]$ such that the product
$[v_iu_{i+2}]\cdot[v_ju_{j+2}u_{j+3}]$ is non-zero. This describes
the multiplicative structure in the cohomology of $\zk$. In
particular, its Betti vector is $(1,0,0,5,5,0,0,1)$. There is also
a similar description for $m>5$, see~\cite[Ex.~7.22]{bu-pa02}.
\end{exa}

Now we are going to compare the bigrading in $H^*(\zk)$ determined
by Theorem~\ref{cohom1} with that coming from the bigraded
cellular structure.

Define the map
\begin{gather*}
  j\colon\k[K]\otimes\Lambda[u_1,\dots,u_m] \to C^*(\zk;\k),\\
  v_i \mapsto D_i^*,\quad
  u_i \mapsto T_i^*.
\end{gather*}
Hence, $j$ maps a monomial $v_\s u_{\tau}$ with
$\s\cap\tau=\varnothing$ to the cellular cochain $\mathcal
T(\s,\tau)^*$ with $\mathcal T(\s,\tau)_D=\s$, \ $\mathcal
T(\s,\tau)_T=\tau$ (note that $\s\in K$) and all other monomials
to zero. The above map respects the bigrading in both algebras,
and one directly checks that it commutes with the differentials.

\begin{thm}\label{quism}
The map $j$ is a quasiisomorphism, that is, induces an isomorphism
$\Tor_{\k[v_1,\dots,v_m]}(\k[K],\k)\to H^*(\zk;\k)$ in the
cohomology.
\end{thm}
\begin{proof}
The map $j$ has an obvious additive right inverse
$i\colon C^*(\zk;\k)\to\k[K]\otimes\L[u_1,\dots,u_m]$ sending
$T(\s,\tau)$ to $v_\s u_{\tau}$. The standard cochain homotopy
operator $s$ for the Koszul resolution establishes a cochain
homotopy equivalence between $\mathrm{id}$ and $ij$, that is,
$ds+sd=\mathrm{id}-ij$. See~\cite[Ch.~7]{bu-pa02} for more
details.
\end{proof}

\begin{cor}\label{altop}
The algebraic bigraded Betti numbers of $\k[K]$ coincide with the
topological bigraded Betti numbers of $\zk$:
$$
  \beta^{-i,2j}(\k[K])=b^{-i,2j}(\zk;\k),\quad i,j=1,\dots,m.
$$
\end{cor}

Now we can summarise the results of Proposition~\ref{frmap},
Lemmas~\ref{tamap} and~\ref{mamap}, Corollary~\ref{eqczk} and
Theorem~\ref{quism} in a statement describing the functorial
properties of the correspondence $K\mapsto\zk$. Let us introduce
the following functors:
\begin{itemize}
  \item $\mathcal Z$, the covariant functor $K\mapsto\zk$ from the
  category of finite simplicial complexes and simplicial maps to the
  category of toric spaces and equivariant maps (the \emph{moment-angle complex
  functor});
  \item $\k[\cdot]$, the contravariant functor $K\mapsto\k[K]$ from
  simplicial complexes to graded $\k$-algebras (the \emph{Stanley--Reisner functor});
  \item $\mbox{\rm Tor-alg}$, the contravariant functor
  $$
    K\mapsto\Tor_{\k[v_1,\dots,v_m]}\bigl(\k[K],\k\bigr)
  $$
  from complexes to bigraded $\k$-algebras (the \emph{$\Tor$-algebra
  functor});
  \item $H^*_T$, the contravariant functor $X\mapsto H^*_T(X;\k)$ from
  toric spaces and equivariant maps to $\k$-algebras
  (the \emph{equivariant cohomology functor});
  \item $H^*$, the contravariant functor $X\mapsto H^*(X;\k)$ from
  spaces to $\k$-algebras (the \emph{ordinary cohomology functor}).
\end{itemize}

\begin{prop} We have the following identities of functors:
$$
  H^*_T\circ\mathcal Z=\k[\cdot],\qquad
  H^*\circ\mathcal Z=\mbox{\rm Tor-alg}.
$$
\end{prop}
The later identity means that for every simplicial map $\phi\colon K_1\to K_2$
the cohomology map $\phi^*_{\mathcal Z}\colon H^*(\mathcal Z_{K_2};\k)\to
H^*(\mathcal Z_{K_1};\k)$ coincides with the induced homomorphism
$\phi_t^*$~{\rm(\ref{toralgmap})} of $\Tor$-algebras.

Also, now we are ready to give a proof of Theorem ~\ref{hoch}.
\begin{proof}[Proof of the Hochster theorem]
First we observe that the cellular cochains of $\zk$ decompose as
\begin{equation}\label{odeco}
  C^{*,\,*}(\zk)=\bigoplus_{\omega\subset[m]}C^{*,\,2\omega}(\zk)
\end{equation}
(as bigraded differential modules), where $C^{*,\,2\omega}(\zk)$
is the subcomplex generated by the cochains $\mathcal
T(\s,\omega\setminus\s)^*$ with $\s\subset\omega$ (remember that
$\mathcal T(\s,\omega\setminus\s)_D=\s$ and $\mathcal
T(\s,\omega\setminus\s)_T=\omega\setminus\s$). It follows that
\begin{equation}\label{bomeg}
  b^{-i,\,2j}(\zk)=
  \sum_{\omega\subset[m]\colon|\omega|=j}b^{-i,\,2\omega}(\zk),
\end{equation}
where $b^{-i,\,2\omega}(\zk):=\dim H^{-i}[C^{*,\,2\omega}(\zk)]$.

The assignment $\s^*\mapsto\mathcal
T(\s,\omega\setminus\s)^*\subset C^{|\s|-|\omega|,\,2\omega}(\zk)$
defines an isomorphism of cochain complexes
$$
  C^*(K_\omega)\to C^{*+1-|\omega|,\,2\omega}(\zk),
$$
where the former is the simplicial cochains of $K_\omega$.
Hence,
$$
  b^{-i,\,2\omega}(\zk)= b^{|\omega|-i-1}(K_\omega).
$$
This together with \eqref{bomeg} and Corollary~\ref{altop} implies
the Hochster theorem.
\end{proof}

\subsection{Generalised Dehn--Sommerville equations}
Here we consider the bigraded Poincar\'e duality for $\zk$. As a
corollary, we deduce linear relations for the number of faces in a
triangulated manifold, generalising the Dehn--Sommerville
equations.

The decomposition \eqref{odeco} gives rise to the coarser
decomposition
$$
  C^{*,\,*}(\zk)=\bigoplus_{p=0}^m C^{*,\,2p}(\zk).
$$
Let us consider the corresponding Euler characteristics:
$$
  \chi_p(\zk):=\sum_{q=0}^m(-1)^q\dim C^{-q,\,2p}(\zk)
  =\sum_{q=0}^m(-1)^qb^{-q,\,2p}(\zk)
$$
and define the generating polynomial $\chi(\zk;t)$ by
$$
  \chi(\zk;t)=\sum_{p=0}^m\chi_p(\zk)t^{2p}.
$$
It turns out that this polynomial can be expressed in terms of the
number of faces of~$K$. Introduce the \emph{$h$-polynomial} of $K$
as $h(t)=h_0+h_1t+\dots+h_nt^n$, where $(h_0,h_1,\ldots,h_n)$ is
the $h$-vector.

\begin{lem}\label{chihv}
For every $K^{n-1}$ it holds that
$$
  \chi(\zk;t)=(1-t^2)^{m-n}h(t^2).
$$
\end{lem}
\begin{proof}
The component $C^{-q,\,2p}(\zk)$ has basis of cochains $\mathcal
T(\s,\tau)^*$ with $\s\in K$, \ $|\s|=p-q$ and $|\tau|=q$. It
follows that
$$
  \dim C^{-q,\,2p}(\zk)=f_{p-q-1}\bin{m-p+q}q.
$$
The rest is a direct calculation using \eqref{hvector}.
See~\cite[Thm.~7.15]{bu-pa02} for details.
\end{proof}

We may regard the standard torus $T^m\subset\C^m$ as a cellular
subcomplex in the moment-angle complex $\zk$. In the same fashion
as we did before for $\zk$, we may define the bigraded Betti
numbers and characteristic polynomials for the pair $(\zk,T^m)$
and for $\zk\setminus T^m$. The proof of the following statement
uses a similar but rather more complicated argument as that of
Lemma~\ref{chihv}.

\begin{lem}\label{chitm}
For every $K^{n-1}$ it holds that
\begin{align*}
  \chi(\zk,T^m;t) &=(1-t^2)^{m-n}h(t^2)-(1-t^2)^m;\\
  \chi(\zk\setminus T^m;t)
  &=(1-t^2)^{m-n}h(t^2)+(-1)^{n-1}h_n(1-t^2)^m.
\end{align*}
\end{lem}

Assume that $K$ is a triangulation of $S^{n-1}$. Then $\zk$ is an
$(m+n)$-dimensional (closed) manifold.

\begin{exe}
Show that the top cohomology group $H^{m+n}(\zk)$ is generated
by the cohomology class of any monomial $v_\s
u_\tau\in\k[K]\otimes\L[u_1,\dots,u_m]$ of bidegree $(-(m-n),2m)$
satisfying $\s\cap\tau=\varnothing$.
\end{exe}

\begin{cor}
The Poincar\'e duality for $\zk$ respects the bigraded structure
in the cohomology. In particular,
$$
  b^{-q,\,2p}(\zk)=b^{-(m-n)+q,\,2(m-p)}(\zk).
$$
\end{cor}

It follows that
$$
  \chi_p(\zk)=(-1)^{m-n}\chi_{m-p}(\zk)\quad\text{and}\quad
  \chi(\zk;t)=(-1)^{m-n}t^{2m}\chi(\zk;\textstyle\frac1t).
$$

\begin{exe}
Deduce the Dehn--Sommerville equations $h_i=h_{n-i}$ from the
above identity and Lemma~\ref{chihv}.
\end{exe}

Now assume that $K$ is a triangulation of a closed
$(n-1)$-dimensional manifold. In this case $\zk$ may fail to be a
manifold, but $\zk\setminus T^m$ still has a homotopy type of a
manifold with boundary. The relative Poincar\'e duality still
regards the bigradings, so we have
\begin{equation}\label{rbgpd}
  \chi(\zk\setminus T^m;t)=
  (-1)^{m-n}t^{2m}\chi(\zk,T^m;\textstyle\frac1t).
\end{equation}

\begin{thm}[Dehn--Sommerville equations for triangulated manifolds]
The following relations hold for the $h$-vector
$(h_0,h_1,\ldots,h_n)$ of any triangulated manifold~$K^{n-1}$:
$$
  h_{n-i}-h_i=(-1)^i\bigl(\chi(K^{n-1})-\chi(S^{n-1})\bigr)
  {\textstyle\binom ni},\quad i=0,1,\ldots,n.
$$
\end{thm}
\begin{proof}
It follows from Lemma \ref{chitm} and \eqref{rbgpd}. See
\cite[Thm.~7.44]{bu-pa02} for details.
\end{proof}
Note that $\chi(K^{n-1})-\chi(S^{n-1})=(-1)^{n-1}(h_n-1)$, so the
above relations are indeed linear equations for the numbers of
faces.

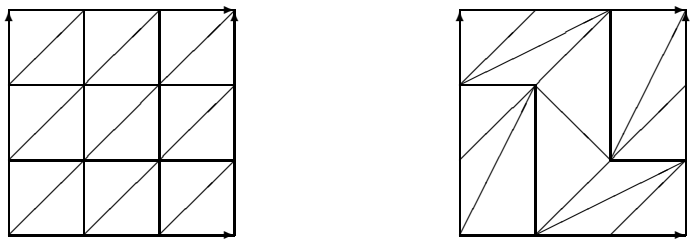
\begin{figure}[h]
  \begin{center}
  \begin{picture}(120,40)(3,0)
  \put(15,10){\vector(1,0){30}}
  \put(15,10){\vector(0,1){30}}
  \put(15,40){\vector(1,0){30}}
  \put(45,10){\vector(0,1){30}}
  \multiput(15,20)(0,10){3}{\line(1,0){30}}
  \multiput(25,10)(10,0){3}{\line(0,1){30}}
  \put(15,30){\line(1,1){10}}
  \put(15,20){\line(1,1){20}}
  \put(15,10){\line(1,1){30}}
  \put(25,10){\line(1,1){20}}
  \put(35,10){\line(1,1){10}}
  \put(0,0){(a)\ $\mb f=(9,27,18)$, $\mb h=(1,6,12,-1)$}
  \put(75,10){\vector(1,0){30}}
  \put(75,10){\vector(0,1){30}}
  \put(75,40){\vector(1,0){30}}
  \put(105,10){\vector(0,1){30}}
  \put(75,10){\line(1,2){10}}
  \put(75,20){\line(1,1){20}}
  \put(75,30){\line(2,1){20}}
  \put(75,30){\line(1,1){10}}
  \put(75,30){\line(1,0){10}}
  \put(85,30){\line(1,-1){10}}
  \put(85,10){\line(0,1){20}}
  \put(85,10){\line(1,1){20}}
  \put(85,10){\line(2,1){20}}
  \put(95,10){\line(1,1){10}}
  \put(95,20){\line(1,0){10}}
  \put(95,20){\line(1,2){10}}
  \put(95,20){\line(0,1){20}}
  \put(67,0){(b)\ $\mb f=(7,21,14)$, $\mb h=(1,4,10,-1)$}
  \end{picture}
  \end{center}
  \caption{``Symmetric" and ``minimal" triangulation of $T^2$}
  \label{tortri}
\end{figure}
\begin{exa}
Consider triangulations of the 2-torus~$T^2$. We have $n=3$, \
$\chi(T^2)=0$. From $\chi(K^{n-1})=1+(-1)^{n-1}h_n$ we deduce
$h_3=-1$. The Dehn--Sommerville equations are
$$
  h_3-h_0=-2,\quad h_2-h_1=6.
$$
For instance, the triangulation on Figure~\ref{tortri}~(a) has
$f_0=9$ vertices, $f_1=27$ edges and $f_2=18$ triangles. The
corresponding $h$-vector is $(1,6,12,-1)$.

On the other hand, a triangulation of $T^2$ with only 7 vertices
can be achieved, see Figure~\ref{tortri}~(b). Note that this
triangulation is neighbourly, i.e. its 1-skeleton is a complete
graph on 7 vertices.
\end{exa}

\begin{exe}
Show that any triangulation of $T^2$ has at least 7 vertices.
\end{exe}

\section{Coordinate subspace arrangements}
Yet another construction of the moment-angle complex $\mathcal
Z_K$ comes from the study of coordinate subspace arrangements.
Coordinate subspace arrangements in $m$-dimensional complex space
are in functorial one-to-one correspondence with simplicial
complexes on $m$ vertices. The complement to a coordinate subspace
arrangement is homotopy equivalent to moment-angle complex
$\mathcal Z_K$ for the corresponding simplicial complex $K$. We
investigate different consequences of this fact. In particular, we
compare the Goresky--MacPherson calculations of Betti numbers of a
coordinate subspace arrangement complement with Hochster's
calculation for the bigraded Betti numbers of the
Tor-algebra.

A \emph{complex coordinate subspace arrangement} is a set
$\mathcal{CA}=\{L_1,\ldots,L_r\}$ of coordinate subspaces in some
$\C^m$. Each coordinate subspace can be written as
\begin{equation}\label{zeroc}
  L_\omega=\{(z_1,\ldots,z_m)\in\C^m\colon z_{i_1}=\cdots=z_{i_k}=0\},
\end{equation}
where $\omega=\{i_1,\ldots,i_k\}$ is a subset of~$[m]$.

For each simplicial complex $K$ on the set $[m]$ we define the
complex coordinate subspace arrangement
$\mathcal{CA}(K):=\{L_\omega\colon\omega\notin K\}$ and
its \emph{complement}
$$
  U(K):=\C^m\setminus\bigcup_{\omega\notin K}L_\omega.
$$
Note that if $L\subset K$ is a subcomplex, then $U(L)\subset
U(K)$. The following observation is straightforward.

\begin{prop}
The assignment $K\mapsto U(K)$ defines a one-to-one
order-preserving correspondence between simplicial
complexes on $[m]$ and coordinate subspace arrangement
complements in $\C^m$.
\end{prop}

\begin{exa}\label{compl}
1. If $K=\D^{m-1}$ then $U(K)=\C^m$.

2. If $K=\partial\D^{m-1}$ then
$U(K)=\C^m\setminus\{0\}$.

3. If $K$ is a disjoint union of $m$ vertices, then $U(K)$ is the
complement in $\C^m$ to the set of all codimension-two coordinate
subspaces $z_i=z_j=0$, \ $1\le i<j\le m$.
\end{exa}

The complement $U(K)$ is invariant under the diagonal action of
$T^m$ on $\C^m$.

\begin{prop}
There is an equivariant deformation retraction $U(K)\to\zk$.
\end{prop}
\begin{proof}
For arbitrary subset $\omega\subset[m]$, set
$$
  U_\omega:=\bigl\{ (z_1,\ldots,z_m)\in \C^m\colon
  |z_i|\ne0\text{ if }i\notin\omega\bigr\}.
$$
Then $U(K)=\bigcup_{\s\in K}U_\s$. Each $B_\sigma$
in~\eqref{bdeco} is a deformation retract of $U_\sigma$, and these
deformation retractions patch together to yield a deformation
retraction $U(K)\to\zk$.
\end{proof}

Hence, we may use the previous results on moment-angle complexes
to calculate the cohomology ring of a coordinate subspace
arrangement complement.

\begin{cor}\label{cohar}
The following isomorphism of graded algebras holds:
$$
  H^{*}\bigl(U(K);\k\bigr)\cong
  \Tor_{\k[v_1,\ldots,v_m]}\bigl(\k[K],\k\bigr)
  =H\bigl[\L[u_1,\ldots,u_m]\otimes\k(K),d\bigr].
$$
\end{cor}

\begin{exa}\label{discr}
Let $K$ and $U(K)$ be as in Example \ref{compl}.3. Then
$\k(K)=\k[v_1,\ldots,v_m]/\mathcal I_K$, where $\mathcal I_K$ is generated by
the monomials $v_iv_j$, \ $i\ne j$.  An easy calculation using
Corollary~\ref{cohar} shows that monomials
$v_{i_1}u_{i_2}u_{i_3}\cdots u_{i_k}$ with $k\ge2$ and $i_p\ne i_q$ for
$p\ne q$ form a basis of cocycles in $\k(K)\otimes\Lambda[u_1,\ldots,u_m]$.
Since $\deg(v_{i_1}u_{i_2}u_{i_3}\cdots u_{i_k})=k+1$, the space of
degree-$(k+1)$ cocycles has dimension $m\binom{m-1}{k-1}$.  The space of
degree-$(k+1)$ coboundaries is $\binom mk$-dimensional and is spanned by
the coboundaries of the form $d(u_{i_1}\cdots u_{i_k})$.  Hence,
\begin{align*}
  &\dim H^{0}\bigl(U(K)\bigr)=1,\quad
  H^{1}\bigl(U(K)\bigr)=H^{2}\bigl(U(K)\bigr)=0,\\
  &\dim H^{k+1}\bigl(U(K)\bigr)=
  m\bin{m-1}{k-1}-\bin mk=(k-1)\bin mk,\quad2\le k\le m,
\end{align*}
and the multiplication in the cohomology is trivial.

In particular, for $m=3$ we have 6 three-dimensional cohomology classes
$[v_iu_j]$, \ $i\ne j$, subject to 3 relations $[v_iu_j]=[v_ju_i]$, and 3
four-dimensional cohomology classes $[v_1u_2u_3]$, $[v_2u_1u_3]$,
$[v_3u_1u_2]$ subject to one relation
$$
  [v_1u_2u_3]-[v_2u_1u_3]+[v_3u_1u_2]=0.
$$
Hence, $\dim H^{3}(U(K))=3$, \ $\dim H^4(U(K))=2$,
and the multiplication is trivial. It can be shown that $U(K)$ in this case
has a homotopy type of a wedge of spheres:
$$
  U(K)\simeq S^3\vee S^3\vee S^3\vee S^4\vee S^4.
$$
\end{exa}

One can define a coordinate subspace as the linear span
of a subset of the standard basis in $\C^m$, rather than by setting some
coordinates to be zero as in~\eqref{zeroc}. This gives an
alternative way to parametrise coordinate subspace arrangements by
simplicial complexes. Namely, we can write
$$
  \mathcal{CA}(K)=\bigl\{ L_\omega\colon\omega\notin K \bigr\}=
  \bigl\{\mathop{\mathrm{span}}\<e_{i_1},\dots,e_{i_k}\>\colon
  \{i_1,\dots,i_k\}\in \widehat{K} \bigr\}
$$
where $\widehat{K}$ is the simplicial complex given by
$$
  \widehat{K}:=\{\omega\subset[m]\colon[m]\setminus\omega\notin K\}.
$$
It is called the \emph{dual complex} of $K$. The cohomology of full
subcomplexes in $K$ is related to the homology of links in $\widehat{K}$
by means of the following combinatorial version of the Alexander duality
theorem.

\begin{thm}[Alexander duality]\label{aldua}
For any simplicial complex $K\ne\D^{m-1}$ on~$[m]$ and simplex
$\s\in\widehat{K}$ it holds that
$$
  \widetilde{H}^i(K_{\widehat\s})\cong
  \widetilde{H}_{m-3-i-|\s|}\bigl(\link_{\widehat{K}}\s\bigr),
$$
where $\widehat{\s}:=[m]\setminus{\s}$. In particular, for
$\sigma=\varnothing$ we get
$$
  \widetilde{H}^i(K)\cong
  \widetilde{H}_{m-3-i}(\widehat{K}),\quad-1\le j\le m-2,
$$
\end{thm}

The Alexander duality and Hochster theorem allow us to give
a simple argument for the \emph{Goresky--MacPherson formula} in
the coordinate subspace arrangement case, providing another
way to calculate the cohomology of the complement. The theorem
below was firstly proved using elaborated algebraic geometry techniques, such
as \emph{stratified Morse theory} and \emph{intersection cohomology}.

\begin{thm}[Goresky--MacPherson]
We have
$$
  \widetilde{H}^i\bigl(U(K)\bigr)=\bigoplus_{\sigma\in\widehat{K}}
  \widetilde{H}_{2m-2|\sigma|-i-2}(\link_{\widehat{K}}\sigma).
$$
\end{thm}
\begin{proof}
Using Corollary~\ref{cohar} to identify $\beta^{-i,2j}(\k[K])$ with
$\dim_\k H^{-i,2j}(U(K))$, we get the following formula from Hochster's
Theorem~\ref{hoch}:
$$
  H^p\bigl(U(K)\bigr)=
  \bigoplus_{\omega\subset[m]}\widetilde{H}^{p-|\omega|-1}(K_\omega).
$$
Non-empty simplices $\omega\in K$ do not contribute to the above sum, since
the corresponding subcomplexes $K_\omega$ are contractible.
Since $\widetilde{H}^{-1}(\varnothing)=\k$, the empty subset of $[m]$ only
contributes $\k$ to $H^0(U(K))$. Hence,
we may rewrite the above formula as
$$
  \widetilde{H}^p\bigl(U(K)\bigr)
  =\bigoplus_{\omega\notin K}\widetilde{H}^{p-|\omega|-1}(K_\omega)
  =\bigoplus_{\sigma\in\widehat{K}}\widetilde{H}^{p-m+|\sigma|-1}
  (K_{\widehat\sigma}).
$$
Using Theorem \ref{aldua}, we identify the latter with
$\oplus_{\sigma\in\widehat{K}}
\widetilde{H}_{2m-2|\sigma|-p-2}(\link_{\widehat{K}}\sigma)$,
as claimed.
\end{proof}

\begin{exe}
Calculate the cohomology of $U(K)$ from Example~\ref{discr}
using the Goresky--MacPherson formula.
\end{exe}


\begin{thebibliography}{9}

\bibitem{br-he98}
Winfried Bruns and J\"urgen Herzog, \emph{Cohen-Macaulay rings},
revised edition, Cambridge Studies in Adv. Math. {\bf 39},
Cambridge Univ. Press, Cambridge, 1998.

\bibitem{bu-pa02}
Victor M. Buchstaber, Taras E. Panov, \emph{Torus actions and
their applications in topology and combinatorics}, University
Lecture Series {\bf 24}, Amer. Math. Soc., Providence, R.I., 2002.

\bibitem{da-ja91}
Michael W. Davis and Tadeusz Januszkiewicz, \emph{Convex
polytopes, Coxeter orbifolds and torus actions}, Duke Math. J.
{\bf 62} (1991), no.~2, 417--451.

\bibitem{ha-ma??}
Akio Hattori and Mikiya Masuda,
\emph{Theory of multi-fans},
preprint (2001), arXiv: math.SG/0106229.

\bibitem{smit67}
Larry Smith, \emph{Homological algebra and the Eilenberg--Moore
spectral sequence}, Trans. Amer. Math. Soc. {\bf 129} (1967),
58--93.

\bibitem{stan96}
Richard P. Stanley, \emph{Combinatorics and Commutative Algebra},
second edition, Progress in Math. {\bf 41}, Birkh\"auser, Boston,
1996.

\end{thebibliography}
\end{document}